\documentclass[12pt]{elsarticle}
\usepackage{amsmath,amsxtra,amssymb,latexsym, amscd,amsfonts,enumerate,ifthen,indentfirst,amstext}
\usepackage{mathrsfs}
\usepackage{longtable}
\usepackage{tablestyles,tabularx, xcolor,multirow}

\usepackage{pgf,tikz,pgfplots}
\pgfplotsset{compat=1.15}
\usepackage{mathrsfs}
\usetikzlibrary{arrows}

\usepackage{smartdiagram}

\usepackage[english]{babel}
\usepackage{amsthm}

\theoremstyle{plain}
\newtheorem{theorem}{\indent\sc Theorem}[section]
\newtheorem{lemma}[theorem]{\indent\sc Lemma}
\newtheorem{corollary}[theorem]{\indent\sc Corollary}
\newtheorem{proposition}[theorem]{\indent\sc Proposition}

\theoremstyle{definition} 
\newtheorem{definition}[theorem]{\indent\sc Definition}
\newtheorem{remark}[theorem]{\indent\sc Remark}

\numberwithin{equation}{section}
\journal{Journal of Mathematical Analysis and Applications }

\begin{document}
	
	\begin{frontmatter}
		\title{ Admissible solutions to augmented 
			nonsymmetric $k-$Hessian type equations I. 
			The $d-$concavity of the $k-$Hessian type functions}  
		
		\author[first]{Tran Van Bang}
		\address[first]{Department of Mathematics,
			Hanoi Pedagogical University No. 2, 
			 Vinh Phuc, Vietnam }
		\author[second]{Ha Tien Ngoan}
	\address[second]{Institute of  Mathematics, Vietnam Academy of Science and Technology, Hanoi, Vietnam }
		\author[third]{Nguyen Huu Tho}
	\address[third]{Department of  Mathematics, Thuyloi University, Hanoi, Vietnam}
		\author[fourth]{Phan Trong Tien}
	\address[fourth]{Department of Mathematics, Quang Binh University,  Quangbinh, Vietnam}
\begin{abstract}
 We establish for $2 \le k \le n-1$ the strict concavity of the function $f_k(\lambda)=\log(\sigma_k(\lambda))$ on a subset of the positive cone $\Gamma_n=\{\lambda=(\lambda_{1}, \lambda_{2}, \cdots,\lambda_{n})\in \mathbb{R}^n; \lambda_j>0,j=1,\cdots, n\}$ 
where $\sigma_{k}(\lambda)$ is the basic symmetric polynomial of degree $k,$ $2 \leq k \leq n.$ Then we apply the result to study the so-called $d-$concavity of the $k-$Hessian type function $F_{k}(R)=\log \left(S_{k}(R)\right),$ where $S_{k}(R)=\sigma_{k}(\lambda(R)), \lambda(R)=	\left(\lambda_{1}, \lambda_{2}, \cdots, \lambda_{n}\right) \in \mathbb{C}^{n}$ is eigenvalue-vector of $R \in \mathbb{R}^{n \times n},$
$R=\omega+\beta, \omega^{T}=\omega, \omega>0, \quad \beta^{T}=-\beta.$ The $d-$concavity 
will be used in our next paper to study the existence of admissible solutions to the
Dirichlet problem for the augmented nonsymmetric
$k-$Hessian type equations.
\end{abstract}

	\end{frontmatter} 
	

\section{Introduction}
It is well-known that by applying the method of continuity, the study of the $k-$Hessian type equations is reduced to get a priory $C^{2,\alpha}-$estimates, that are strongly connected to the concavity of some related $k-$Hession type functions, whose variables are square matrices.

First example of those equations is the symmetric Monge-Ampere type equations (\cite{1}-\cite{9})
$$
\det\left(D^{2} u-A(x, u, Du)\right)=f(x, u, Du),
$$
where $A\left(x, z, p\right)$ is a symmetric matrix. The corresponding $n-$Hessian type function is the function $g_{1}(\omega)=\sqrt[n]{\det \omega}$ or $g_{2}(\omega)=\log(\det \omega),$
that is concave on the set of positive definite matrices $\omega,$ i.e. $\omega^{T}=\omega, \quad \omega > 0.$

Second example is the symmetric $k-$Hessian type equations (\cite{10}-\cite{13})
$$
S_{k}\left(D^{2} {u}-A(x, u, Du)\right)=f(x,u, Du),$$
where $A(x,z,p)$ is a symmetric matrix, 
\begin{equation}\label{1.1}
	S_{k}(\omega)=\sigma_{k}(\lambda(\omega)),
\end{equation}
$\lambda(\omega)=(\lambda_{1}, \lambda_{2}, \cdots, \lambda_{n})$ is eigenvalues vector of $\omega,$
$\sigma_{k}(\lambda)$  is the basic symmetric polynomial of degree $k, 1\le k\le n,$

\begin{equation}\label{1.2}
	 \sigma_k(\lambda)=\sum_{1\le i_1<\cdots<i_k\le n}  \lambda_{i_1}\lambda_{i_2}\cdots\lambda_{i_k}.
\end{equation}
The corresponding $k-$Hessian type function is the function $g(\omega)=\sqrt[k]{S_{k}(\omega)},$ that is concave (\cite{2}) on a connected part of the following set: 
\begin{equation*}\label{1.3}
\Gamma_k:=\{\lambda=(\lambda_{1},\cdots,\lambda_{n})\in \mathbb{R}^n; \sigma_j(\lambda)>0, j=1,2,\cdots,k\}.
\end{equation*}
The third example of $k-$Hessian type equations is the non-symmetric Monge-Ampere type equations (\cite{14},\cite{15}):
$$
\det\left(D^{2} u-A(x, u, D u)-B(x, u, D u)\right)=f(x,u,Du),
$$
where $A(x, z, p)$ is symmetric, $B(x, z, p)$ is skew symmetric. The investigating these equations
is related to the $d-$concavity of the function
$$
F(R)=\log (\det R)=\log \left(S_{n}(R)\right)
$$
in the following sense
\begin{equation}\label{1.4}
F\left(R^{(1)}\right)-F\left(R^{(0)}\right) \leq \sum_{i,j=1}^{n} \frac{\partial F\left(R^{(0)}\right)}{\partial R_{ij}}\left(R_{i j}^{(1)}-R_{i j}^{(0)}\right)+d
\end{equation}
on the set
\begin{equation*}\label{1.5}
		\begin{split}
			D_{\delta,\mu}=\{R=[R_{i j}]_{n\times n}\in \mathbb{R}^{n\times n};&\ R=\omega +\beta,\ \omega ^T=\omega,\ \omega>0,\\
			&\beta^T=-\beta,\ \|\beta\|\leq \mu,\ \mu\leq \delta\lambda_{\min}(\omega )\}.
		\end{split}
\end{equation*}
where $0 \leq \delta<1,  \mu \ge 0, \lambda_{\min }(\omega)$ is the least
eigenvalue of $\omega,$ $\|\beta\|$ is the operator norm of $\beta$
and $d \geqslant 0,$ that depends only on $n, \delta$ and does not depend on $\mu, R^{(0)}, R^{(1)}.$ To prove \eqref{1.4} in \cite{14} it has been established that the following inequalities are true for second-order differential $d^2 F(R,M)$ of $F(R).$
\begin{equation}\label{1.6}
	d^{2} F(R, P) \leq -C_{1}|\tilde{\tilde{P}}|^{2},
\end{equation}
\begin{equation}\label{1.7}
d^{2} F(R, M) \leq  C_{2}|\tilde{\tilde{Q}}|^{2},
\end{equation}
 where $M=P+Q \in \mathbb{R}^{n\times n},   P^T=P,\  Q^{T}=-Q,$
$ R=\omega+\beta \in D_{\delta, \mu}, \omega=C^{-1} D C, C^T=C^{-1}, D=\operatorname{diag}\left(\lambda_{1}, \cdots, \lambda_{n}\right),$
$\lambda=\left(\lambda_{1}, \lambda_{2}, \cdots, \lambda_{n}\right) \in \Gamma_{n},   \widetilde{P}=C P C^{-1},   \tilde{\tilde{P}}=D^{-\frac{1}{2}} \widetilde{P} D^{-\frac{1}{2}},$ $P=[P_{ij}],$ $|P|^2=\sum_{i,j=1}^{n}|P_{ij}|^2; C_1, C_2$ are positive, and do not depend on $R, P, Q$ and $\mu.$ This means that the inequalities \eqref{1.4}, \eqref{1.6} and \eqref{1.7} are uniform with respect to $R^{(0)}, R^{(1)}$ and $R,$ that vary in $D_{\delta,\mu}.$ We note that \eqref{1.6} plays an essential role in getting \eqref{1.7}. We note also that: $\|P\|\le |P|\le \sqrt{n}\|P\|,$ $\|Q\|\le |Q|\le \sqrt{n}\|Q\|.$

The aim of our study here and in our next paper (\cite{16}) is the nonsymmetric augmented $k-$Hessian type equations:
\begin{equation*}\label{1.8}
	S_{k}\left(D^{2} {u}-A(x, u,Du)-B(x, u, Du)\right)=f\left(x, u, D u\right),
\end{equation*}
where $2 \leq  k \leq n, A(x, z, p)$ is symmetric, $B(x, z, p)$ is skew-symmetric. So, in this paper, we would like to study the $d-$concavity of the following corresponding $k-$Hessian type function:
\begin{equation*}\label{1.9}
	F_k(R):=\log\left(S_k(R)\right), \quad  2\le k\le n,
\end{equation*}
on the set $D_{\delta,\mu}.$ We note that, $ S_k(R)\in \mathbb{R},$ and  moreover, if $R\in D_{\delta, \mu}$ then $S_k(R)>0.$ (Corollary \ref{hq2.1}).

Since $S_n(R)=\det R,$ the case $k=n$ had been therefore studied in \cite{14}, \cite{15}, as it was described above. In this paper, for the cases $2\le k\le n-1$ we will prove (Theorem \ref{dl5.3}) the $d-$concavity of the function $F_k(R)$ on some convex and unbounded subset $D_{\delta,\mu,\gamma_k}$ of the set $D_{\delta,\mu}.$

The paper is organized as follows. In Section 2 we try to reduce general case $R\in D_{\delta,\mu}$ to the case $\widetilde{R} =D+\widetilde{\beta}\in D_{\delta,\mu},$ where $D=\operatorname{diag}(\lambda_{1},\lambda_{2},\cdots,\lambda_{n})$ and $D>0.$ In Section 3 we continue to reduce this case to the case $\widetilde{R}=D$ and the increment $\widetilde{P}$ is also diagonal. In Section 4 we establish the strict concavity of the usual function $f_k(\lambda)=\log\left(\sigma_{k}(\lambda)\right)$ on a subset $\Sigma_{(\gamma_k)}$ of the positive cone $\Gamma_n.$ Then we show that the second-order differential $d^2f_k(\lambda,\xi)$ of $f_k(\lambda)$ where $2\le k\le n-1$ satisfies the following inequality
\begin{equation}\label{1.10}
d^{2} f_{k}(\lambda, \xi) \leq -C_{3}|\widetilde{\xi}|^{2}, \quad  \xi \in \mathbb{R}^{n}	
\end{equation}
where $\widetilde{\xi}=\left( \frac{\xi_1}{\lambda_1},\cdots,\frac{\xi_n}{\lambda_n} \right), C_3>0$ and does not depend on $\lambda,$ that varies in a convex and unbounded subset $\Sigma_{(\gamma_k)}$ of $\Gamma_{n}.$ In the last Section 5 we apply \eqref{1.10} to  establish the corresponding inequalities \eqref{1.6}, \eqref{1.7} for the function $F_{k}(R)$ when $R\in D_{\delta,\mu,\gamma_k},$ where $\delta>0$ is sufficiently small, depends on $n,k,\gamma_k$ and does not depend on $\mu.$ The $d-$concavity of $F_k(R)$ is a consequence of the corresponding inequality \eqref{1.7}. 
\section{The representation of second-order differential}
For the matrix $R=[R_{ij}]_{n\times n}$ and  the indices $i_1i_2\cdots i_k,$ which are always assumed such that, $1\le i_1<i_2<\cdots <i_k\le n,$ we introduce following notations: 
\begin{equation}\label{2.1}
R_{i_1\cdots i_k}=[R_{i_pi_q}]_{p,q=1}^k,	
\end{equation}
\begin{equation}\label{2.2}
G_{i_1\cdots i_k}(R)=\det(R_{i_1\cdots i_k}),	
\end{equation}
\begin{equation}\label{2.3}
H_{i_1\cdots i_k}(R)=\log(G_{i_1\cdots i_k}(R)).
\end{equation}
We know (\cite{13}) that
\begin{equation}\label{2.4}
S_{k}(R)=\sigma_{k}(\lambda(R))=\sum_{i_1,\cdots, i_k} G_{i_1 \cdots  i_{k}}(R)	,
\end{equation}
where $\displaystyle\sum_{i_1,\cdots, i_k} $ stands for $\displaystyle\sum_{1\le i_1<\cdots< i_k\le n}. $ \\
Suppose $R=\omega+\beta, \omega=[\omega_{ij}]_{n\times n}, \beta=[\beta_{ij}]_{n\times n},$ $\omega^T=\omega, \omega>0, \beta^T=-\beta.$
\begin{proposition}[\cite{14}, Proposition 3]\label{md2.1}
Suppose $R=\omega+\beta\in D_{\delta,\mu}.$ Then the following assertions hold for any index $i_1i_2\cdots i_k:$

\begin{itemize}
	\item[(i)]  $\omega^T_{i_1\cdots i_k}=\omega_{i_1\cdots i_k},\quad  \omega_{i_1\cdots i_k}>0$ and
				$\lambda_{\min}(\omega_{i_1\cdots i_k})	\ge \lambda_{\min}(\omega);$
		\item[(ii)]  $\beta^{T}_{i_1\cdots i_k}=-\beta_{i_1\cdots i_k}, \quad \|\beta_{i_1\cdots i_k}\|\le \mu,  \quad \det \beta_{i_1\cdots i_k}\ge 0;$
	\item[(iii)]  $R_{i_1\cdots i_k}=\omega_{i_1\cdots i_k} +\beta_{i_1\cdots i_k}\in D_{\delta,\mu};$
		\item[(iv)]  $\det(R_{i_1\cdots i_k})\ge \det \omega_{i_1\cdots i_k}+
		\det \beta_{i_1 \cdots i_k}\ge \det \omega_{i_1\cdots i_k}>0.$
\end{itemize}
\end{proposition}
\begin{corollary}\label{hq2.1}
 If $R\in D_{\delta,\mu}$ then $S_k(R)>0,\quad k=1,2,\cdots,n.$
\end{corollary}
\begin{proposition}[\cite{14}, Proposition 7] \label{md2.2}
	Suppose $R=\left[R_{ij}\right]_{n\times n}\in D_{\delta,\mu}.$ Set 
$$H(R)=\log( \det R),\quad R^{-1}=[(R^{-1})_{ij}]_{n\times n}.$$
Then the following assertion holds 
\begin{equation}\label{2.6}
	\frac{\partial H(R)}{\partial R_{ij} }=(R^{-1})_{ji}.
\end{equation}
\end{proposition}

From \eqref{2.2}, \eqref{2.3} and \eqref{2.6}   we have 
\begin{equation}\label{2.8}
	 \dfrac{\partial H_{i_1\cdots i_k}(R)}{\partial R_{i j}}=\sum_{p,q=1}^k (R_{i_1\cdots i_k})_{i_q i_p}^{-1}\delta_{i i_p}\delta_{j i_q};
\end{equation}
\begin{equation}\label{2.9}
	\dfrac{\partial^2 H_{i_1\cdots i_k}(R)}{\partial R_{i j}\partial R_{\ell m}}=\sum_{p,q,r,s=1}^k \frac{\partial^2 H_{i_1\cdots i_k}(R)}{\partial R_{i_p i_q}\partial R_{i_r i_s}}\delta_{i i_p}\delta_{j i_q} \delta_{\ell i_r}\delta_{m i_s}.
\end{equation}
It follows from \eqref{2.4} and \eqref{2.3} that 
$$F_k(R)=\log (S_k(R))=\log \left(\sum_{i_1,\cdots,i_k} e^{H_{i_1\cdots i_k}(R)}\right).$$
The derivatives of $F_k(R)$ are given by
\begin{equation*}\label{2.10}
	\frac{\partial F_k(R)}{\partial R_{ij}}=\frac{1}{S_k(R)}\sum_{i_1,\cdots,i_k} G_{i_1\cdots i_k}(R)  \dfrac{\partial H_{i_1\cdots i_k}(R)}{\partial R_{i j}},
\end{equation*}
and 
 \begin{equation}\label{2.11}
 	 \begin{split}
	\dfrac{\partial^2 F_k(R)}{\partial R_{i j}\partial R_{\ell m}}=&-\frac{1}{S_k^2(R)}\underset{j_1,\cdots,j_k}{\sum_{i_1,\cdots,i_k}} G_{i_1\cdots i_k}(R)G_{j_1\cdots j_k}(R)\dfrac{\partial H_{i_1\cdots i_k}(R)}{\partial R_{i j}}\dfrac{\partial H_{j_1\cdots j_k}(R)}{\partial R_{\ell m}},\\
	&+\frac{1}{S_k(R)}\sum_{i_1,\cdots,i_k} G_{i_1\cdots i_k}(R)\dfrac{\partial H_{i_1\cdots i_k}(R)}{\partial R_{i j}}\dfrac{\partial H_{i_1\cdots i_k}(R)}{\partial R_{\ell m}}\\
	&+\frac{1}{S_k(R)}\sum_{i_1,\cdots,i_k} G_{i_1\cdots i_k}(R)\dfrac{\partial^2 H_{i_1\cdots i_k}(R)}{\partial R_{i j}\partial R_{\ell m}}.
\end{split}\end{equation}
Suppose $M=[M_{ij}]_{n\times n}\in \mathbb{R}^{n\times n}.$ From \eqref{2.11} we have the following representation for second-order differential of $F_k(R)$ at $R\in D_{\delta,\mu}:$
\begin{equation}\label{2.12}
	\begin{split}
	d^2F_k(R,M):=&\sum_{i,j,\ell,m=1}^n \dfrac{\partial^2 F_k(R)}{\partial R_{i j}\partial R_{\ell m}}M_{i j}M_{\ell m}\\
	=&- \frac{1}{S_k^2(R)}\Big[\underset{i, j=1}{\sum_{i_1,\cdots,i_k}^n}G_{i_1\cdots i_k}(R)\dfrac{\partial H_{i_1\cdots i_k}(R)}{\partial R_{i j}}M_{i j}\Big]^2\\
	&+\frac{1}{S_k(R)}\sum_{i_1,\cdots,i_k} G_{i_1\cdots i_k}(R)\Big[\sum_{i,j=1}^n\dfrac{\partial H_{i_1\cdots i_k}(R)}{\partial R_{i j}}M_{i j}\Big]^2\\
	&+\frac{1}{S_k(R)}\sum_{i_1,\cdots,i_k} G_{i_1\cdots i_k}(R)\sum_{i,j,\ell,m=1}^n\dfrac{\partial^2 H_{i_1\cdots i_k}(R)}{\partial R_{i j}\partial R_{\ell m}}M_{i j}M_{\ell m}.
\end{split}\end{equation}
Since $d^{2} F_{k}(R, M)$ is a symmetric quadratic form with respect to $M,$ then it is easy to  get the following proposition
\begin{proposition}\label{md2.3}
		Suppose $M=P+Q,$ where $P^T=P, Q^T=-Q.$ Then we have 
		\begin{equation}\label{2.13}
			d^2F_k(R,M)=d^2 F_k(R,P)+d^2 F_k(R,Q)+2 H_k(R,P,Q),
		\end{equation}	
		where
		\begin{equation*}\label{2.14}
			H_{k}(R, P, Q)=\sum_{i, j,\ell,m =1}^n \dfrac{\partial^{2} F_{k}(R)}{\partial R_{ij} \partial R_{\ell m}}  P_{i j} Q_{\ell m}.
		\end{equation*}
	\end{proposition}
Now we reduce the general case $R\in D_{\delta,\mu}$ to a simpler one, where the symmetric part of $R$ is a diagonal matrix. 
\begin{proposition}\label{md2.4}
Suppose $R=\omega+\beta\in D_{\delta,\mu}$ and $\omega=C^{-1}DC,$ where $C$ is an orthogonal matrix and $D={\rm diag}(\lambda_{1},\lambda_{2},\cdots,\lambda_{n})>0.$ Set 
\begin{equation}\label{2.15}
	\widetilde{R}=CR C^{-1}=D+\widetilde{\beta}, \quad\tilde{\beta}=C \beta C^{-1},
\end{equation}
\begin{equation}\label{2.16}
\widetilde{M}=CMC^{-1}.
\end{equation}
Then we have 
\begin{equation}\label{2.17}
	d^{2} F_{k}(R, M)=d^{2} F_{k}(\widetilde{R}, \widetilde{M}).
\end{equation}
\end{proposition}
\begin{proof}
	
	Suppose $C=[C_{ij}]_{n\times n}$ and $C^{-1}=C^T=\left[(C^{-1})_{ij}\right]_{n\times n}.$ From \eqref{2.15} and \eqref{2.16}  we have 
	\begin{equation}\label{2.18}
		\widetilde{R}_{pq}=\sum_{i,j=1}^{n} C_{pi} R_{ij}(C^{-1})_{jq}, \quad p,q=1,\cdots,n,
	\end{equation}
\begin{equation*}\label{2.19}
	\widetilde{M}_{p q}=\sum_{i, j=1}^{n} C_{p i} M_{i j}  (C^{-1})_{j q},\quad  p, q=1, \cdots, n.
\end{equation*}
From \eqref{2.18} it follows 
$$\frac{\partial F_k(R)}{\partial R_{ij}}=\sum_{p,q=1}^{n} C_{pi} \frac{\partial F_k(\widetilde{R})}{\partial \widetilde{R}_{pq}}(C^{-1})_{jq},$$
$$\frac{ \partial^2 F_k(R)}{ \partial R_{ij}\partial R_{\ell m}}=\sum_{p,q=1}^{n} C_{pi} \left(\sum_{r,s=1}^n C_{r\ell} \frac{\partial^2 F_k(\widetilde{R}) }{ \partial\widetilde{R}_{pq} \partial\widetilde{R}_{rs} }(C^{-1})_{ms}\right)(C^{-1})_{jq}.$$
Therefore 
\begin{equation*}
	\begin{split}
	 \sum_{i,j,\ell, m=1}^{n}  \dfrac{\partial^{2} F_{k}(R)}{\partial R_{i j} \partial R_{\ell m}}  M_{i j} M_{\ell m}&  =\sum_{p,q,r,s=1}^{n}  \dfrac{\partial^{2} F_{k}(\widetilde{R})}{\partial \widetilde{R}_{pq}\partial \widetilde{R}_{rs}} \left( \sum_{i,j =1}^{n} C_{pi} M_{ij} (C^{-1})_{jq} \right)\left( \sum_{\ell, m =1}^{n} C_{r\ell} M_{\ell m} (C^{-1})_{ms}\right)\\
	&=   \sum_{p,q,r,s=1}^{n}  \dfrac{\partial^{2} F_{k}(\widetilde{R})}{\partial \widetilde{R}_{pq}\partial \widetilde{R}_{rs}} \widetilde{M}_{pq}\widetilde{M}_{rs}=d^2 F_k(\widetilde{R},\widetilde{M}).
\end{split}\end{equation*}
	\end{proof}
\section{Estimates for second-order differential}
In this Section we try to reduce the case $\widetilde{R}=D+\widetilde{\beta}$ to the case where $\widetilde{R}=D={\rm diag}\left(\lambda_{1}, \lambda_{2}, \cdots, \lambda_{n}\right)>0$ and the symmetric part $\widetilde{P}$ of the increment $\widetilde{M}=\widetilde{P}+\widetilde{Q}$ is also a diagonal matrix. We note that if $R\in D_{\delta,\mu}$ then $\widetilde{R}=D+\widetilde{\beta}\in D_{\delta,\mu}.$
\subsection{Some preliminary lemmas}
First of all we prove some lemmas. We suppose always that $\widetilde{R}=D+\widetilde{\beta}\in D_{\delta,\mu}.$
\begin{lemma}[\cite{14}, Proposition 6]\label{bd3.1}
	Suppose 
	$$\widetilde{\sigma}=D^{-\frac{1}{2}}\widetilde{\beta} D^{-\frac{1}{2}},$$
	where $\widetilde{\beta}$ is defined by \eqref{2.15}. 
	Then for any indices $i_1\cdots i_k$ we have
\begin{itemize}
	\item[(i)] 
	\begin{equation}\label{3.1}
		\widetilde{\sigma}_{i_1\cdots i_k}= D^{-\frac{1}{2}}_{i_1\cdots i_k} \widetilde{\beta}_{i_1\cdots i_k} D^{-\frac{1}{2}}_{i_1\cdots i_k};
	\end{equation} 
	\item[(ii)]  \begin{equation}\label{3.2}
		\|	\widetilde{\sigma}_{i_1\cdots i_k} \|\le \delta,
	\end{equation}
where $\widetilde{\sigma}_{i_1\cdots i_k}$ is defined by \eqref{2.1}. Moreover, 
there exists $K_{i_1\cdots i_k}(\widetilde{R})\in \mathbb{R}^{k\times k}$ such that
	\item[(iii)] 
		\begin{equation}\label{3.3}
\left|K_{i_{1} \cdots i_k}  (\widetilde{R})\right| \leq \frac{\sqrt{k} \delta^{2}}{\left(1-\delta^{2}\right)}; (E_{i_1\cdots i_k}+\widetilde{\sigma}_{i_1\cdots i_k})^{-1}=E_{i_1\cdots i_k}+ K_{i_1\cdots i_k}
\end{equation}
\item[(iv)]  
\begin{equation}\label{3.4}
\frac{ (\widetilde{R}_{i_1\cdots i_k})^{-1}+\left[ (\widetilde{R}_{i_1\cdots i_k})^{-1}\right]^T}{2}=D^{-\frac{1}{2}}_{i_1\cdots i_k} \left( E_{i_1\cdots i_k}+K_{i_1\cdots i_k}\right) D^{-\frac{1}{2}}_{i_1\cdots i_k};
\end{equation}	
\item[(v)] \begin{equation}\label{3.5}
	\frac{ (\widetilde{R}_{i_1\cdots i_k})^{-1}-\left[ (\widetilde{R}_{i_1\cdots i_k})^{-1}\right]^T}{2}=D^{-\frac{1}{2}}_{i_1\cdots i_k}\left(-\widetilde{\sigma}_{i_1\cdots i_k} \right) \left( E_{i_1\cdots i_k}+K_{i_1\cdots i_k}\right) D^{-\frac{1}{2}}_{i_1\cdots i_k},
\end{equation}	
where $E_{i_{1} \cdots i_k}$ is the unit matrix in $\mathbb{R}^{k \times k}.$\end{itemize}
	\end{lemma}
\begin{proof}
We prove only the inequality $\left\|\widetilde{\sigma}_{i_1\cdots i_k}\right\| \leq\delta.$ Indeed, from the equality in \eqref{3.1} and Proposition \ref{md2.1} we have 
$$\left\|\widetilde{\sigma}_{i_1\cdots i_k}\right\|\le  \left\|\widetilde{\beta}_{i_1\cdots i_k}\right\|\left\|D^{-\frac{1}{2}}_{i_1\cdots i_k}\right\|^2\le \frac{\mu}{\lambda_{\min}(D_{i_1\cdots i_k})}\le \frac{\mu}{\lambda_{\min}(\omega)}\le \delta.$$
\end{proof}
\begin{lemma}\label{bd3.2}
	For any indices $i_1\cdots i_k$ there exists a function $h_{i_1\cdots i_k}(\widetilde{R})$ such that 
	\begin{itemize}
		\item[(i)] 
		\begin{equation}\label{3.6}
			0\le h_{i_1\cdots i_k}(\widetilde{R}) \le \left(2^{[\frac{k}{2}]}-1 \right)\delta^2,
		\end{equation}
		\item[(ii)] 
		\begin{equation}\label{3.7}
			G_{i_1\cdots i_k}(\widetilde{R})= \left(1+ h_{i_1\cdots i_k}(\widetilde{R}) \right)\det (D_{i_1\cdots i_k}),
		\end{equation}	
	where $	G_{i_1\cdots i_k}(\widetilde{R})$ is defined by \eqref{2.2}.
	\end{itemize}
\end{lemma}
\begin{proof}
	Suppose $\widetilde{\sigma}_{i_1\cdots i_k}$ is defined by \eqref{3.1}. It is a skew-symmetric matrix of order $k.$ We denote its eigenvalues by $i{\eta_1}, i{\eta_2},\cdots, i{\eta_k},$ where $\eta_j\in \mathbb{R},$ $\eta_2=-\eta_1,\cdots, \eta_{2[\frac{k}{2}]}=-\eta_{2[\frac{k}{2}]-1}, \eta_k=0$ if $k$ is odd and $|\eta_j|\le \delta.$ Then, from the proof of Proposition 3 in \cite{14}, it follows that
	$$G_{i_1\cdots i_k}(\widetilde{R})=\left(1+\eta^2_1\right)\left(1+\eta^2_3\right)\cdots \left(1+\eta^2_{2[\frac{k}{2}]-1}\right) \det(D_{i_1\cdots i_k}).$$
	We choose 
	\begin{equation}\label{3.8}
		h_{i_1\cdots i_k}(\widetilde{R})=\left(1+\eta^2_1\right)\left(1+\eta^2_3\right)\cdots \left(1+\eta^2_{2[\frac{k}{2}]-1}\right) -1\ge0.
	\end{equation}
Then we have immediately \eqref{3.7}. Since $|\eta_j|\le \delta<1$ it follows that
\begin{equation*}
	h_{i_1\cdots i_k}(\widetilde{R})=\sum_{j=1}^{[\frac{k}{2}]}	\sigma_j\left(\eta^2_1,\eta^2_3,\cdots, \eta^2_{2[\frac{k}{2}]-1} \right)\le \sum_{j=1}^{[\frac{k}{2}]} \binom{[\frac{k}{2}]}{j}\delta^2=\left(2^{[\frac{k}{2}]}-1\right)\delta^2.
\end{equation*}
	\end{proof}
From Lemma \ref{bd3.2} following corollaries follow.
\begin{corollary}
	We set 
	\begin{equation*}\label{3.9}
		h_k(\widetilde{R}) =\frac{1}{S_k(D)} \sum_{i_1,\cdots,i_k} h_{i_1\cdots i_k}(\widetilde{R}) \det (D_{i_1\cdots i_k}).
	\end{equation*}
Then the following assertions are true 
	\end{corollary}
 	\begin{itemize}
 	\item[(i)] 
 	\begin{equation}\label{3.10}
 		0\le h_k (\widetilde{R})\le \left(2^{[\frac{k}{2}]}-1 \right)\delta^2,
 			\end{equation}
 			\item[(ii)] 
 			\begin{equation}\label{3.11}
 				S_k (\widetilde{R})=\left(1+h_k(\widetilde{R})\right) S_k(D).
 			\end{equation}
 \end{itemize}
\begin{corollary}
For any indices $i_{1} \cdots i_k$ there exists
a function $g_{i_{1}\cdots i_{k}}(\widetilde{R})$ such that
	\begin{itemize}
	\item[(i)] 
\begin{equation}\label{3.12}
\left|g_{i_{1}\cdots i_{k}}(\widetilde{R})\right| \leq \left(2^{\left[\frac{k}{2}\right]}-1\right) \delta^{2},
\end{equation}
\item[(ii)]
\begin{equation}\label{3.13}
	\frac{G_{i_{1}\cdots i_{k}}(\widetilde{R})  }{S_{k}(\widetilde{R})}=\frac{G_{i_{1}\cdots i_{k}}(D)}{S_{k}(D)}\left(1+g_{i_{1}\cdots i_{k}}(\widetilde{R})\right).
\end{equation}
\end{itemize}
	\end{corollary}
\begin{proof}
It is sufficient to choose 
$$g_{i_{1}\cdots i_{k}}(\widetilde{R}) =\frac{1+h_{i_{1}\cdots i_{k}}(\widetilde{R})}{1+h_k(\widetilde{R})}-1=\frac{h_{i_{1}\cdots i_{k}}(\widetilde{R})-h_k(\widetilde{R})}{1+h_k(\widetilde{R})}.
$$
Then \eqref{3.12} follows from \eqref{3.8} and \eqref{3.10}.
	\end{proof}
\begin{lemma}\label{bd3.3}
Suppose $\widetilde{P}^T=\widetilde{P}$ and $K_{i_{1}\cdots i_{k}}	(\widetilde{R})$ is the matrix in Lemma \ref{bd3.1}. Then 
	\begin{itemize}
	\item[(i)] 
	\begin{equation}\label{3.14}
		{\rm Tr}\left[ \left(\widetilde{R}_{i_{1}\cdots i_{k}}\right)^{-1}\widetilde{P}_{i_{1}\cdots i_{k}}\right] ={\rm Tr} \left(\tilde{\tilde{P}}_{i_{1}\cdots i_{k}}\right)+{\rm Tr} \left(K_{i_{1}\cdots i_{k}}\tilde{\tilde{P}}_{i_{1}\cdots i_{k}}\right),
	\end{equation}
\item[(ii)]  
\begin{equation}\label{3.15}
	\left|{\rm Tr}\left(\tilde{\tilde{P}}_{i_{1}\cdots i_{k}}\right)\right| \le \sqrt{k}| \tilde{\tilde{P}}_{i_{1}\cdots i_{k}}| \le \sqrt{k}| \tilde{\tilde{P}} |,\end{equation}
\item[(iii)] 
\begin{equation}\label{3.16}
	\left|{\rm Tr} \left(K_{i_{1}\cdots i_{k}}\tilde{\tilde{P}}_{i_{1}\cdots i_{k}}\right)\right|\le \left| K_{i_{1}\cdots i_{k}}\right| \left| \tilde{\tilde{P}}_{i_{1}\cdots i_{k}}\right| \le \frac{\sqrt{k}\delta^2}{(1-\delta^2)}|\tilde{\tilde{P}}|,
\end{equation}
 where ${\rm Tr}(M)=\sum_{j=1}^n M_{jj}$ for $M=[M_{ij}]_{n\times n},$ $\tilde{\tilde{P}}= D^{-\frac{1}{2}}\tilde{P} D^{-\frac{1}{2}}.$
\end{itemize}
	\end{lemma}
\begin{proof}
	Since $\tilde{P}_{i_{1}\cdots i_{k}}$ is symmetric and  $\tilde{\tilde{P}}_{i_{1}\cdots i_{k}}=D^{-\frac{1}{2}}_{i_{1}\cdots i_{k}} \widetilde{P}_{i_{1}\cdots i_{k}}D^{-\frac{1}{2}}_{i_{1}\cdots i_{k}},$ then 
	$${\rm Tr}\left[(\widetilde{R}_{i_{1}\cdots i_{k}})^{-1}  \widetilde{P}_{i_{1}\cdots i_{k}} \right]= {\rm Tr}\left[\frac{(\widetilde{R}_{i_{1}\cdots i_{k}})^{-1} +[(\widetilde{R}_{i_{1}\cdots i_{k}})^{-1}]^T}{2}  \widetilde{P}_{i_{1}\cdots i_{k}} \right] $$
	then \eqref{3.14} follows from \eqref{3.4}. 
	
	Suppose $\tilde{\tilde{P}}_{i_{1}\cdots i_{k}}=C^{-1}_{i_{1}\cdots i_{k}}\widetilde{D}_{i_{1}\cdots i_{k}}C_{i_{1}\cdots i_{k}},$ where $C_{i_{1}\cdots i_{k}}$ is an orthogonal matrix, $\widetilde{D}_{i_{1}\cdots i_{k}}={\rm diag}(\mu_1,\cdots,\mu_k), \mu_i\in \mathbb{R}.$ We have 
	$${\rm Tr}\left(\tilde{\tilde{P}}_{i_{1}\cdots i_{k}} \right)={\rm Tr}(\widetilde{D}_{i_{1}\cdots i_{k}}) =\sum_{j=1}^k \mu_j,$$
	$$\left|{\rm Tr}\left(\tilde{\tilde{P}}_{i_{1}\cdots i_{k}} \right)\right|^2\le k \sum_{j=1}^k \mu_j^2=k\left|\tilde{\tilde{P}}_{i_{1}\cdots i_{k}} \right|^2,$$
	then  \eqref{3.15} follows. The inequality \eqref{3.16} follows from \eqref{3.3}. Here we have used the fact that $\left|\tilde{\tilde{P}}_{i_{1}\cdots i_{k}} \right|\le \left|\tilde{\tilde{P}} \right|.$
	\end{proof}
\begin{lemma}
	Suppose $\widetilde{Q}^T=-\widetilde{Q}.$ Then 
	\begin{align}
		&(i)\ \ \	{\rm Tr}\left[ \left(\widetilde{R}_{i_{1}\cdots i_{k}} \right)^{-1} \widetilde{Q}_{i_{1}\cdots i_{k}} \right]={\rm Tr} \left[ -\widetilde{\sigma}_{i_1\cdots i_k}\tilde{\tilde{Q}}_{i_{1}\cdots i_{k}} \right]+ {\rm Tr} \left[ -\widetilde{\sigma}_{i_1\cdots i_k} K_{i_{1}\cdots i_{k}}\tilde{\tilde{Q}}_{i_{1}\cdots i_{k}} \right], \label{3.17}\\
&	(ii) \ \ 
	\left| {\rm Tr} \left[ -\widetilde{\sigma}_{i_1\cdots i_k}\tilde{\tilde{Q}}_{i_{1}\cdots i_{k}} \right] \right| \le \left| \widetilde{\sigma}_{i_1\cdots i_k}\right| \left|\tilde{\tilde{Q}}_{i_{1}\cdots i_{k}}| \right| \le \sqrt{k}\delta |\tilde{\tilde{Q}}|,\label{3.18}\\
&(iii)\ \left|  {\rm Tr} \left(-\widetilde{\sigma}_{i_1\cdots i_k} K_{i_{1}\cdots i_{k}}\tilde{\tilde{Q}}_{i_{1}\cdots i_{k}} \right)  \right|\le \left| \widetilde{\sigma}_{i_1\cdots i_k}\right| |K_{i_{1}\cdots i_{k}}| \left|\tilde{\tilde{Q}}_{i_{1}\cdots i_{k}} \right|  
		\le \sqrt{k}\delta \frac{\sqrt{k}\delta^2 }{(1-\delta^2)} |\tilde{\tilde{Q}}|. \label{3.19}
\end{align} 
\end{lemma}
\begin{proof}
	It follows from \eqref{3.2} that
	\begin{equation}\label{3.20T}
		\left| \widetilde{\sigma}_{i_1\cdots i_k}  \right|\le \sqrt{k} 	\left\| \widetilde{\sigma}_{i_1\cdots i_k}  \right\| \le \sqrt{k}  \delta.
	\end{equation}
The proof of the lemma is analogous to that of Lemma \ref{bd3.3}. The inequalities \eqref{3.18}, 	\eqref{3.19} follow from \eqref{3.20T}. Here we have used $\left|\tilde{\tilde{Q}}_{i_{1}\cdots i_{k}}\right|\le |\tilde{\tilde{Q}}|.$
	\end{proof}
\begin{lemma}[\cite{14}, Propositions  13 and 14] \label{bd3.5}
		Suppose $\widetilde{\sigma}_{i_1\cdots i_k}$ is defined by \eqref{3.2}, $i\eta_1, i\eta_2, \cdots, i\eta_k$ are 
		 its eigenvalues, where $ \eta_{j} \in \mathbb{R}, \eta_{2}=-\eta_{1}, \cdots, \eta_{2 [\frac{k}{2}]}
		  =-\eta_{2 [\frac{k}{2} ]-1}, \eta_{k}=0 $ if $k$ is odd.

Suppose $\widetilde{\sigma}_{i_1\cdots i_k}=C^{*}_{i_1\cdots i_k} D^{(k)} C_{i_1\cdots i_k},$ where $C_{i_1\cdots i_k}$ is an unitary matrix, $D^{(k)}$ is a diagonal, $D^{(k)}=(i\eta_1, i\eta_2, \cdots, i\eta_k), $ $\tilde{\tilde{\tilde{P}}}_{i_1\cdots i_k}=C_{i_1\cdots i_k}\tilde{\tilde{P}}_{i_1\cdots i_k}C^{*}_{i_1\cdots i_k},$ $\tilde{\tilde{\tilde{Q}}}_{i_1\cdots i_k}=C_{i_1\cdots i_k}\tilde{\tilde{Q}}_{i_1\cdots i_k}C^{*}_{i_1\cdots i_k}.$ Then for $\widetilde{P}^T=\widetilde{P}, \widetilde{Q}^T=-\widetilde{Q}$ we have 
\begin{align}
&(i)\ \ 
	\sum_{p,q,r,s=1}^{k} \frac{\partial^2 H_{i_1\cdots i_k}(\widetilde{R} )}{\partial\widetilde{R}_{i_p i_q} \partial\widetilde{R}_{i_r i_s}
	 }	 \widetilde{P}_{i_p i_q}  \widetilde{P}_{i_r i_s}
 =-\sum_{p,q=1}^{k} \frac{(1-\eta_p \eta_q)}{(1+\eta^2_p)(1+\eta^2_q)}\left|\tilde{\tilde{\tilde{P}}}_{i_pi_q
} \right|^2,& & \label{3.20}\\
&(ii)\ 	
	\sum_{p,q,r,s=1}^{k} \frac{\partial^2 H_{i_1\cdots i_k}(\widetilde{R} )}{\partial\widetilde{R}_{i_p i_q} \partial\widetilde{R}_{i_r i_s}
		}	 \widetilde{Q}_{i_p i_q}  \widetilde{Q}_{i_r i_s}
		=\sum_{p,q=1}^{k} \frac{(1-\eta_p \eta_q)}{(1+\eta^2_p)(1+\eta^2_q)}\left|\tilde{\tilde{\tilde{Q}}}_{i_pi_q
		} \right|^2.&  &\label{3.21}
	\end{align} 
	\end{lemma}
From Lemma \ref{bd3.5} we get the following
\begin{lemma}
	Let all conditions of Lemma \ref{bd3.5} fulfil. Then the following assertions hold 	
\begin{align}
	&(i)\ \ 
		\sum_{p,q,r,s=1}^{k} \frac{\partial^2 H_{i_1\cdots i_k}(\widetilde{R} )}{\partial\widetilde{R}_{i_p i_q} \partial\widetilde{R}_{i_r i_s}
		}	 \widetilde{P}_{i_p i_q}  \widetilde{P}_{i_r i_s}
		\le - \left|\tilde{\tilde{P}}_{i_1\cdots i_k}\right|^2+4\delta^2 \left|\tilde{\tilde{P}}_{i_1\cdots i_k}\right|^2,& & \label{3.22}\\
&(ii)\ 	
		\sum_{p,q,r,s=1}^{k} \frac{\partial^2 H_{i_1\cdots i_k}(\widetilde{R} )}{\partial\widetilde{R}_{i_p i_q} \partial\widetilde{R}_{i_r i_s}
		}	 \widetilde{Q}_{i_p i_q}  \widetilde{Q}_{i_r i_s}
\le \left|\tilde{\tilde{Q}}_{i_1 \cdots i_k}\right|^2+4\delta^2 \left|\tilde{\tilde{Q}}_{i_1 \cdots i_k}\right|^2.& & \label{3.23}
\end{align}	
	\end{lemma}
\begin{proof}
Let 
\begin{equation}\label{3.24}
	m_{pq}(\widetilde{R}) =\frac{1-\eta_p\eta_q}{(1+\eta^2_p)(1+\eta^2_q)}-1=-\frac{\eta_p\eta_q+\eta^2_{p}+\eta^2_q+\eta^2_{p}\eta^2_q}{(1+\eta^2_p)(1+\eta^2_q)}.
\end{equation}	
Since $|\eta_j|\le \delta$ then 
\begin{equation}\label{3.25}
	\left| m_{pq} (\widetilde{R})\right| \le 3\delta^2+\delta^4\le 4\delta^2.
\end{equation}	
Then \eqref{3.22} follows from \eqref{3.20}, \eqref{3.24}, \eqref{3.25}  and $|\tilde{\tilde{\tilde{P}}}_{i_1 \cdots i_k}|=|\tilde{\tilde{P}}_{i_1 \cdots i_k}|.$ The inequality \eqref{3.23} is proved by the same way.
	\end{proof}
\subsection{Some estimates}
\begin{proposition}\label{md3.1}
	Assume that $\widetilde{P}^T=\widetilde{P}.$ Then 
	\begin{equation}\label{3.26}
		d^2 F_k(\widetilde{R},\widetilde{P})\le 	d^2 F_k(D,\widetilde{P})+C_4 \delta^2|\tilde{\tilde{P}}|^2.
	\end{equation}
	\end{proposition}
\begin{proof}
	From \eqref{2.12}, \eqref{2.8} and \eqref{2.9}, it follows
	\begin{equation}\label{3.27}
		\begin{split}
			d^2 F_k(\widetilde{R},\widetilde{P})=& -\frac{	1
			}{S^2_k(\widetilde{R})}\left[ \sum_{i_1,\cdots,i_k} G_{i_1\cdots i_k} (\widetilde{R}) {\rm Tr} \left[( \widetilde{R}_{i_1\cdots i_k})^{-1} \widetilde{P}_{i_1\cdots i_k}\right]\right]^2\\
		&+\frac{	1
		}{S_k(\widetilde{R})} \sum_{i_1,\cdots,i_k} G_{i_1\cdots i_k} (\widetilde{R}) \left[ {\rm Tr} \left( \widetilde{R}_{i_1\cdots i_k}\right)^{-1} \widetilde{P}_{i_1\cdots i_k}\right]^2\\
	&+\frac{	1
	}{S_k(\widetilde{R})} \sum_{i_1,\cdots,i_k} G_{i_1\cdots i_k} (\widetilde{R}) \sum_{p,q,r,s=1}^{k} \frac{\partial^2 H_{i_1\cdots i_k}(\widetilde{R})}{\partial \widetilde{R}_{i_pi_q}\partial \widetilde{R}_{i_ri_s}}\widetilde{P}_{i_pi_q}\widetilde{P}_{i_ri_s}.
		\end{split}
	\end{equation}
From \eqref{3.14} we have 
\begin{equation*}\label{3.28}
		\begin{split}
	\sum_{i_1,\cdots, i_k} G_{i_1\cdots i_k} (\widetilde{R}).& {\rm Tr} \left[ (\tilde{R}_{i_1\cdots i_k})^{-1}  \tilde{P}_{i_1 \cdots i_k} \right]\\
&=	\sum_{i_1,\cdots, i_k} G_{i_1\cdots i_k} (\widetilde{R}). {\rm Tr} \left( \tilde{\tilde{P}}_{i_1\cdots i_k}\right) + \sum_{i_1,\cdots, i_k} G_{i_1\cdots i_k} (\widetilde{R}). {\rm Tr} \left[ K_{i_1\cdots i_k}\tilde{\tilde{P}}_{i_1\cdots i_k}\right],
		\end{split}
\end{equation*}
\begin{equation}\label{3.29}
	\begin{split}
		\Big[\sum_{i_1,\cdots, i_k} G_{i_1\cdots i_k}& (\widetilde{R}) {\rm Tr} \left[ (\tilde{R}_{i_1\cdots i_k})^{-1}  \tilde{P}_{i_1\cdots i_k} \right] \Big]^2 \\
		=&	\left[\sum_{i_1,\cdots, i_k} G_{i_1\cdots i_k} (\widetilde{R}) {\rm Tr} \left( \tilde{\tilde{P}}_{i_1\cdots i_k}\right)\right]^2 + \left[\sum_{i_1,\cdots, i_k} G_{i_1\cdots i_k} (\widetilde{R}). {\rm Tr} \left[ K_{i_1\cdots i_k}\tilde{\tilde{P}}_{i_1\cdots i_k}\right]\right]^2\\
		&+2 \left[\sum_{i_1,\cdots, i_k} G_{i_1\cdots i_k} (\widetilde{R}) {\rm Tr} \left( \tilde{\tilde{P}}_{i_1\cdots i_k}\right)\right]  \left[\sum_{i_1,\cdots, i_k} G_{i_1\cdots i_k} (\widetilde{R}). {\rm Tr} \left[ K_{i_1\cdots i_k}\tilde{\tilde{P}}_{i_1\cdots i_k}\right]\right].
	\end{split}
\end{equation}	
From \eqref{3.13}, \eqref{3.12}, \eqref{3.15}, \eqref{3.16} it follows that
\begin{equation}\label{3.30}
	\begin{split}
	-\frac{	1
	}{S^2_k(\widetilde{R})}\Big[ \sum_{i_1,\cdots, i_k} &G_{i_1\cdots i_k} (\widetilde{R}) {\rm Tr} \left[(\tilde{\tilde{P}}_{i_1\cdots i_k})\right]\Big]^2\\
& \le - \frac{	1
}{S^2_k(D)}\left[ \sum_{i_1,\cdots, i_k} G_{i_1\cdots i_k} (D) {\rm Tr} (\tilde{\tilde{P}}_{i_1\cdots i_k})\right]^2 
+(2^{[\frac{k}{2}]}-1)\delta^2 |\tilde{\tilde{P}}|^2,
	\end{split}
\end{equation}
\begin{equation}\label{3.31}
	\begin{split}
	\left|\frac{	1
		}{S_k(\widetilde{R})}\left[ \sum_{i_1,\cdots, i_k} G_{i_1 \cdots i_k} (\widetilde{R}) {\rm Tr} \left[K_{i_1 \cdots i_k}\tilde{\tilde{P}}_{i_1 \cdots i_k}\right]\right]\right|	& \le \sum_{i_1,\cdots, i_k}\frac{G_{i_1 \cdots i_k} (\widetilde{R})}{S_k(\widetilde{R})} |K_{i_1 \cdots i_k} |\tilde{\tilde{P}}_{i_1 \cdots i_k}|\\
	&\le \sum_{i_1,\cdots, i_k}\frac{G_{i_1 \cdots i_k} (\widetilde{R})}{S_k(\widetilde{R})} \frac{\sqrt{k}\delta^2}{(1-\delta^2)}|\tilde{\tilde{P}}|\\
	&=\frac{\sqrt{k}\delta^4}{(1-\delta^2)}|\tilde{\tilde{P}}|;
	\end{split}
\end{equation}
\begin{equation}\label{3.32}
	\begin{split}
	\left| 	\frac{	1
		}{S_k(\widetilde{R})}\left[ \sum_{i_1,\cdots, i_k} G_{i_1 \cdots i_k} (\widetilde{R}) {\rm Tr} \left[(\tilde{\tilde{P}}_{i_1 \cdots i_k})\right]\right]\right| & \le \sum_{i_1,\cdots, i_k}\frac{G_{i_1 \cdots i_k} (\widetilde{R})}{S_k(\widetilde{R})} 	\left| {\rm Tr} \left[(\tilde{\tilde{P}}_{i_1 \cdots i_k})\right] \right|= \sqrt{k}  |\tilde{\tilde{P}}|.
	\end{split}
\end{equation}
Here we have used equalities:
$$\sum_{i_1,\cdots, i_k}\frac{G_{i_1 \cdots i_k} (\widetilde{R})}{S_k(\widetilde{R})}  =\sum_{i_1,\cdots, i_k}\frac{G_{i_1 \cdots i_k} (D)}{S_k(D)} =1.$$
From \eqref{3.29}-\eqref{3.32}  it follows that 
\begin{equation}\label{3.33}
	\begin{split}
	  	-\frac{	1
		}{S^2_k(\widetilde{R})}\Big[ \sum_{i_1,\cdots, i_k} G_{i_1 \cdots i_k} (\widetilde{R})& {\rm Tr} \left[( \tilde{R}_{i_1 \cdots i_k})^{-1} \tilde{P}_{i_1 \cdots i_k}\right]\Big]^2 \\
	&  \le 	-\frac{	1
	}{S^2_k(D)}    \Big[\sum_{i_1,\cdots, i_k}   G_{i_1 \cdots i_k} (D)    {\rm Tr} (\tilde{\tilde{P}}_{i_1 \cdots i_k}) \Big]^2 +C_6 \delta^2  |\tilde{\tilde{P}}|^2,
	\end{split}
\end{equation}
where $$C_6=(2^{[\frac{k}{2}]}-1)+\frac{k \delta^2}{(1-\delta^2)^2}+\frac{2k}{1-\delta^2}.$$
 Now we estimate the second term on the right-hand side of \eqref{3.27}. From \eqref{3.14} we have 
\begin{equation}\label{3.34}
	\begin{split}
\big[{\rm Tr}&\big[(\widetilde{R}_{i_1 \cdots i_k})^{-1}  \widetilde{P}_{i_1 \cdots i_k} \big]\big]^{2}\\
&=	 \left[{\rm Tr} ( \tilde{\tilde{P}}_{i_1 \cdots i_k})\right]^2
	+ 	\left[{\rm Tr} (K_{i_1 \cdots i_k} \tilde{\tilde{P}}_{i_1 \cdots i_k})\right]^2+2 \left[{\rm Tr} ( \tilde{\tilde{P}}_{i_1 \cdots i_k})\right] \left[{\rm Tr} (K_{i_1 \cdots i_k} \tilde{\tilde{P}}_{i_1 \cdots i_k})\right].
 	\end{split}
\end{equation}
From \eqref{3.13}, \eqref{3.12}, \eqref{3.15}, \eqref{3.16} it follows that 
\begin{equation}\label{3.35}
	\begin{split}
	 \frac{	1}{S_k(\widetilde{R})}  \sum_{i_1,\cdots, i_k} G_{i_1 \cdots i_k} (\widetilde{R})& \left[{\rm Tr}  (\tilde{\tilde{P}}_{i_1 \cdots i_k})\right]^2\\
	 &\le  \frac{	1}{S_k(D)}  \sum_{i_1,\cdots, i_k} G_{i_1 \cdots i_k} (D) \left[{\rm Tr}  (\tilde{\tilde{P}}_{i_1 \cdots i_k})\right]^2+(2^{[\frac{k}{2}]}-1)\delta^2|\tilde{\tilde{P}}|^2,
	\end{split}
\end{equation}
\begin{equation}\label{3.36}
	\begin{split}
	\Big|\frac{	1}{S_k(\widetilde{R})}  	\sum_{i_1,\cdots, i_k} G_{i_1 \cdots i_k} (\widetilde{R})& \left[{\rm Tr}  (K_{i_1 \cdots i_k}\tilde{\tilde{P}}_{i_1 \cdots i_k})\right]^2 \Big|\\
		&\le \sum_{i_1,\cdots, i_k}\frac{G_{i_1 \cdots i_k} (\widetilde{R})}{S_k(\widetilde{R})}  
		\frac{k\delta^4}{(1-\delta^2)^2} |\tilde{\tilde{P}}|^2=\frac{k\delta^4}{(1-\delta^2)^2} |\tilde{\tilde{P}}|^2,
	\end{split}
\end{equation}
\begin{equation}\label{3.37}
	\begin{split}
		\Big|\frac{1}{S_k(\widetilde{R})}  \sum_{i_1,\cdots, i_k} G_{i_1 \cdots i_k} (\widetilde{R}) &\left[{\rm Tr}(\tilde{\tilde{P}}_{i_1 \cdots i_k})\right]
		\left[{\rm Tr}  (K_{i_1 \cdots i_k}\tilde{\tilde{P}}_{i_1 \cdots i_k})\right] \Big|\\
		& \le \sum_{i_1,\cdots, i_k}\frac{G_{i_1 \cdots i_k} (\widetilde{R})}{S_k(\widetilde{R})} \sqrt{k} |\tilde{\tilde{P}}| \frac{\sqrt{k}\delta^2}{(1-\delta^2)}|\tilde{\tilde{P}}|=\frac{k\delta^2}{(1-\delta^2)}|\tilde{\tilde{P}}|^2.
	\end{split}
\end{equation}
From \eqref{3.34}-\eqref{3.37} we conclude that 
\begin{equation}\label{3.38}
	\begin{split}
		 \frac{1}{S_k(\widetilde{R})}  \sum_{i_1,\cdots, i_k} G_{i_1 \cdots i_k} (\widetilde{R}) &\left[{\rm Tr}(\tilde{R}_{i_1 \cdots i_k})^{-1} \tilde{P}_{i_1 \cdots i_k}\right]^2\\
		 & \le \frac{1}{S_k(D)}  \sum_{i_1,\cdots, i_k} G_{i_1 \cdots i_k} (D) \left[{\rm Tr}(\tilde{\tilde{P}}_{i_1 \cdots i_k})\right]^2+C_7\delta^2|\tilde{\tilde{P}}|^2,
	\end{split}
\end{equation}
where $$C_7=(2^{[\frac{k}{2}]}-1)+\frac{k\delta^2}{(1+\delta^2)^2}+\frac{2k}{1-\delta^2}=C_6.$$
Now we estimate the last term on the right-hand side of \eqref{3.27}.

Since $G_{i_1\cdots i_k} (\widetilde{R})>0$ and $S_k(\widetilde{R})>0,$ then from \eqref{3.22} we get 
\begin{equation}\label{3.39}
	\begin{split}
	 \frac{1}{S_k(\widetilde{R})} 	& \sum_{i_1,\cdots, i_k} G_{i_1 \cdots i_k} (\widetilde{R})
		\sum_{p,q,r,s=1}^{k} \frac{\partial^2 H_{i_1\cdots i_k}(\widetilde{R})}{\partial \widetilde{R}_{i_p  i_q}\partial \widetilde{R}_{i_ri_s}}\widetilde{P}_{i_pi_q}\widetilde{P}_{i_ri_s}\\
				& \le -\frac{1}{S_k(\widetilde{R})}  \sum_{i_1,\cdots, i_k} G_{i_1 \cdots i_k} (\widetilde{R}) |\tilde{\tilde{P}}_{i_1 \cdots i_k}|^2+4\delta^2|\tilde{\tilde{P}}|^2\\
				&=-\sum_{i_1,\cdots, i_k}\frac{G_{i_1 \cdots i_k} (D)}{S_k(D)}(1+g_{i_1 \cdots i_k}(\widetilde{R})) |\tilde{\tilde{P}}_{i_1 \cdots i_k}|^2+4\delta^2|\tilde{\tilde{P}}|^2\\
				&\le -\sum_{i_1,\cdots, i_k}\frac{G_{i_1 \cdots i_k} (D)}{S_k(D)}  |\tilde{\tilde{P}}_{i_1 \cdots i_k}|^2+4\delta^2|\tilde{\tilde{P}}|^2+4\delta^2|\tilde{\tilde{P}}|^2\\
				&=-\frac{1}{S_k(D)}    \sum_{i_1,\cdots, i_k}G_{i_1 \cdots i_k} (D)	\sum_{p,q,r,s=1}^{k} \frac{\partial^2 H_{i_1 \cdots i_k}(D)}{\partial \widetilde{R}_{i_pi_q}\partial \widetilde{R}_{i_ri_s}}\widetilde{P}_{i_pi_q}\widetilde{P}_{i_ri_s} +8\delta^2|\tilde{\tilde{P}}|^2.
	\end{split}
\end{equation}
So, from \eqref{3.33}, \eqref{3.38} and \eqref{3.39},  \eqref{3.26} follows with $C_4=2C_6+8.$
	\end{proof}
\begin{proposition}\label{md3.2} 
Suppose $ \widetilde{Q}^{T}=-\widetilde{Q}.$  Then 
\begin{equation}\label{3.40}
	d^{2} F_{k}(\widetilde{R}, \widetilde{Q}) \leq  C_{8}|\tilde{\tilde{Q}}|^{2},
\end{equation}
where  $C_{8}=1+C_{4} \delta.$
\end{proposition}
\begin{proof}
	The proof is analogous to that of Proposition \ref{md3.1}. In this case we use \eqref{3.17}, \eqref{3.18}, \eqref{3.19}, \eqref{3.7}, \eqref{3.13},  \eqref{3.21} and \eqref{3.23}.	
	\end{proof}
\begin{proposition}\label{md3.3}
	 Suppose $ \widetilde{P}^{T}=\widetilde{P}, \widetilde{Q}^{T}=-\widetilde{Q}.$ Then 
	 \begin{equation}\label{3.41}
	 	H_k(\widetilde{R}, \widetilde{P}, \widetilde{Q})=\sum_{i,j,\ell, m=1}^{k} \frac{\partial^2 F_k(\widetilde{R})}{\partial \widetilde{R}_{ij}\partial \widetilde{R}_{\ell m}}\widetilde{P}_{ij}\widetilde{Q}_{\ell m} \le C_9\delta |	 	\tilde{\tilde{P}}| |\tilde{\tilde{Q}}|,
	 \end{equation}
 where $$C_9=2k\left[\frac{1+(\sqrt{k}-1)\delta^2}{(1-\delta^2)^2}+ \frac{1}{1+\delta^2}\right].$$
\end{proposition}
\begin{proof}
	From \eqref{3.41}, \eqref{3.11} and \eqref{3.6} it follows that 
	 \begin{equation*}
		\begin{split}
	H_k(\widetilde{R}, \widetilde{P}, \widetilde{Q})=&-\left[\sum_{i_1,\cdots, i_k}\frac{G_{i_1 \cdots i_k} (\widetilde{R})}{S_k(\widetilde{R})}	{\rm Tr}((\widetilde{R}_{i_1 \cdots i_k})^{-1}\widetilde{P}_{i_1 \cdots i_k}) \right]\left[\sum_{j_1, \cdots, j_k}\frac{G_{j_1 \cdots j_k} (\widetilde{R})}{S_k(\widetilde{R})}	{\rm Tr}((\widetilde{R}_{j_1 \cdots j_k})^{-1}\widetilde{Q}_{j_1 \cdots j_k}) \right]  \\
	&+ \sum_{i_1,\cdots, i_k}\frac{G_{i_1 \cdots i_k} (\widetilde{R})}{S_k(\widetilde{R})}	\left[{\rm Tr}((\widetilde{R}_{i_1 \cdots i_k})^{-1}\widetilde{P}_{i_1 \cdots i_k})  \right] \left[  {\rm Tr}((\widetilde{R}_{i_1 \cdots i_k})^{-1}\widetilde{Q}_{i_1 \cdots i_k}) \right]\\
	&+ \sum_{i_1,\cdots, i_k}\frac{G_{i_1 \cdots i_k} (\widetilde{R})}{S_k(\widetilde{R})}	
		\sum_{p,q,r,s=1}^{k} \frac{\partial^2 H_{i_1 \cdots i_k}(\widetilde{R})}{\partial \widetilde{R}_{i_pi_q}\partial \widetilde{R}_{i_ri_s}}\widetilde{P}_{i_pi_q}\widetilde{Q}_{i_ri_s}. 
		\end{split}
	\end{equation*}
	Then we have 
 \begin{equation}\label{3.42}
	\begin{split}
			H_k(\widetilde{R}, \widetilde{P}, \widetilde{Q}) \le& \left[\sum_{i_1,\cdots, i_k}\frac{G_{i_1 \cdots i_k} (\widetilde{R})}{S_k(\widetilde{R})} \left|
			{\rm Tr}((\widetilde{R}_{i_1 \cdots i_k})^{-1}\widetilde{P}_{i_1 \cdots i_k})\right|\right] \left[\sum_{j_1, \cdots, j_k}\frac{G_{j_1 \cdots j_k} (\widetilde{R})}{S_k(\widetilde{R})} \left|
			{\rm Tr}((\widetilde{R}_{j_1 \cdots j_k})^{-1}\widetilde{Q}_{j_1 \cdots j_k})\right|\right] \\
		&+ \sum_{i_1,\cdots, i_k}\frac{G_{i_1 \cdots i_k} (\widetilde{R})}{S_k(\widetilde{R})}	\left|{\rm Tr}((\widetilde{R}_{i_1 \cdots i_k})^{-1}\widetilde{P}_{i_1 \cdots i_k})  \right| \left|  {\rm Tr}((\widetilde{R}_{i_1 \cdots i_k})^{-1}\widetilde{Q}_{i_1 \cdots i_k}) \right|\\
		&+ \sum_{i_1,\cdots, i_k}\frac{G_{i_1 \cdots i_k} (\widetilde{R})}{S_k(\widetilde{R})}	
		\sum_{p,q,r,s=1}^{k} \frac{\partial^2 H_{i_1 \cdots i_k}(\widetilde{R})}{\partial \widetilde{R}_{i_pi_q}\partial \widetilde{R}_{i_ri_s}}\widetilde{P}_{i_pi_q}\widetilde{Q}_{i_ri_s}. 
	\end{split}
\end{equation}
 We will estimate each term on the right-hand side of \eqref{3.42}. From \eqref{3.14}-\eqref{3.16} we obtain 
 \begin{equation*}\label{3.43}
 	\left|{\rm Tr}\left(\left(\widetilde{R}_{i_1 \cdots i_k}\right)^{-1} \widetilde{P}_{i_1 \cdots i_k}\right)\right|\le \left[\sqrt{k}+\frac{\sqrt{k} \delta^{2}}{\left(1-\delta^{2}\right)}\right] |\tilde{\tilde{P}}|,
 \end{equation*}
and from \eqref{3.17}-\eqref{3.19} we have
\begin{equation*}\label{3.44}
		\left|{\rm Tr}\left(\left(\widetilde{R}_{i_1 \cdots i_k}\right)^{-1} \widetilde{Q}_{i_1 \cdots i_k}\right)\right|\le \left[\sqrt{k}\delta+\frac{\sqrt{k} \delta^{3}}{\left(1-\delta^{2}\right)}\right] |\tilde{\tilde{Q}}|.
\end{equation*}
Since $$\sum_{i_1,\cdots, i_k}\frac{G_{i_1\cdots i_k} (\widetilde{R})}{S_k(\widetilde{R})}=1,$$
 both first and second terms  of the right-hand side of \eqref{3.42}  are not greater than 
 $$k\left(1+\frac{\delta^{2}}{1-\delta^{2}}\right)\left(1+\frac{\sqrt{k} \delta^{2}}{1-\delta^{2}}\right) \delta|\tilde{\tilde{P}}| \, |\tilde{\tilde{Q}}|.$$	

 To estimate the third term we need the following 
 \begin{lemma}[\cite{14}, Proposition 15]
 	Suppose $\widetilde{R}=D+\widetilde{\beta}=[\widetilde{R}_{ij}]_{n\times n}\in D_{\delta,\mu}.$ If $\widetilde{P}^T=\widetilde{P}, \widetilde{Q}^T=-\widetilde{Q}, $ $H(\widetilde{R})=\log (\det \widetilde{R}),$ then 
 	\begin{equation}\label{3.45}
 			\sum_{i,j,\ell, m} \frac{\partial^2 H(\widetilde{R})}{\partial \widetilde{R}_{ij}\partial \widetilde{R}_{\ell m}}\widetilde{P}_{ij}\widetilde{Q}_{\ell m}  			\le \frac{2n\delta}{1+\delta^2} |\tilde{\tilde{P}}| |\tilde{\tilde{Q}}|,
 	\end{equation}
 	where 
 	$$\tilde{\tilde{P}}=D^{-\frac{1}{2}} \widetilde{P}D^{-\frac{1}{2}}, \tilde{\tilde{Q}}=D^{-\frac{1}{2}} \widetilde{Q}D^{-\frac{1}{2}}.$$ 
 	\end{lemma}
 We apply \eqref{3.45}   to each function $ H_{i_1 \cdots i_k}(\widetilde{R})=\log (\det \widetilde{R}_{i_1 \cdots i_k})$ and obtain 
  \begin{equation}\label{3.46}
  		\sum_{p,q,r,s=1}^{k} \frac{\partial^2 H_{i_1 \cdots i_k}(\widetilde{R})}{\partial \widetilde{R}_{i_pi_q}\partial \widetilde{R}_{i_ri_s}}\widetilde{P}_{i_pi_q}\widetilde{Q}_{i_ri_s} \le \frac{2k\delta}{1+\delta^2} |\tilde{\tilde{P}}_{i_1 \cdots i_k}| |\tilde{\tilde{Q}}_{i_1 \cdots i_k}| \le \frac{2k\delta}{1+\delta^2} |\tilde{\tilde{P}} | |\tilde{\tilde{Q}} |.
  \end{equation}
From \eqref{3.46}  it follows that the third term on the right-hand side of \eqref{3.42} is not greater than $\frac{2k\delta}{1+\delta^2} |\tilde{\tilde{P}} | |\tilde{\tilde{Q}} |.$ This completes proof of Proposition \ref{md3.3}.
\end{proof}
 Now we will estimate $ d^{2} F_{k}(\widetilde{R}, \widetilde{P}+\widetilde{Q}).$   From 
\eqref{2.13}, \eqref{3.26}, \eqref{3.40} and \eqref{3.41} we get the following.
\begin{proposition}\label{md3.4}
	Suppose $\widetilde{R}=D+\widetilde{\beta}\in D_{\delta,\mu},$ $P^T=P,$ $Q^T=-Q.$ Then for $\widetilde{M}=\widetilde{P}+\widetilde{Q}$ we have
\begin{equation*}\label{3.47}
d^{2} F_{k}(\widetilde{R}, \widetilde{M}) \leq d^{2} F_{k}(D, \widetilde{P})+\left(C_{4}+1\right) \delta^{2}|\tilde{\tilde{P}}|^{2}+\left(C_{9}^{2}+C_{8}\right)|\tilde{\tilde{Q}}|^{2},
\end{equation*}
 where $C_4, C_8, C_9$ depend only on $k, \delta,$ do not depend on $\mu,$ and they  are bounded when $\delta$ tends to zero.
	\end{proposition}
\subsection{Reducing $\widetilde{P}$ to the case of being diagonal}
We now reduce the study of $d^2F_k(D,\widetilde{P})$ to the case where $\widetilde{P}$ is a diagonal matrix. So, suppose 
$D={\rm diag}\left(\lambda_{1}, \lambda_{2}, \cdots, \lambda_{n}\right),$ $\widetilde{P}=[\widetilde{P}_{ij}]_{n\times n},$ $\widetilde{P}^T=\widetilde{P},$ $\tilde{\tilde{P}}=D^{-\frac{1}{2}} \widetilde{P}D^{-\frac{1}{2}}.$ We set 
$$ \widetilde{U} ={\rm diag}\left(\widetilde{P}_{11}, \widetilde{P}_{22}, \cdots, \widetilde{P}_{nn}\right)={\rm diag}\left(\widetilde{U}_{11}, \widetilde{U}_{22}, \cdots, \widetilde{U}_{nn}\right)$$
whose diagonal coincides with that of $\widetilde{P}.$

\begin{proposition}\label{md3.5}
	The following equality is true
\begin{equation}\label{3.48}
	d^{2} F_{k}(D,\widetilde{P}) = d^{2} F_{k}(D, \widetilde{U})-\frac{1}{S_k(D)}\sum_{i_1,\cdots,i_k}G_{i_1\cdots i_k}(D) \sum_{p,q=1\atop p\ne q}^{k}| \tilde{\tilde{P}}_{i_p i_q}|^2.
\end{equation}	
	\end{proposition}
\begin{proof}
	From \eqref{3.27} and \eqref{3.20} we have 
 \begin{equation}\label{3.49}
	\begin{split}
	 	d^{2} F_{k}(D,\widetilde{P})=	& -\frac{1}{\sigma_k^2(\lambda)}\left[ \sum_{i_1, \cdots, i_k} \lambda_{i_1}\cdots \lambda_{i_k} \left(\sum_{p=1}^{k} \frac{\tilde{P}_{i_p i_p}}{\lambda_{i_p}}  \right)\right]^2\\
		&+\frac{1}{\sigma_k(\lambda)}\sum_{i_1, \cdots, i_k}\lambda_{i_1}\cdots \lambda_{i_k}  \left[ \sum_{p=1}^{k} \frac{\tilde{P}_{i_p i_p}}{\lambda_{i_p}}  \right]^2\\
		&-\frac{1}{\sigma_k(\lambda)}  \sum_{i_1,\cdots,i_k}G_{i_1\cdots i_k}(D) \sum_{p,q=1 }^{k}| \tilde{\tilde{P}}_{i_p i_q}|^2,
			\end{split}
\end{equation}
from which we obtain \eqref{3.48}.
	\end{proof}
We write now the explicit form of $	d^{2} F_{k}(D,\widetilde{U}).$
From \eqref{3.49} we have 
\begin{equation}\label{3.50}
	\begin{split}
	  	d^{2} F_{k}(D,\widetilde{U})=	&-\frac{1}{\sigma_k^2(\lambda)}\left[ \sum_{i_1,\cdots, i_k} \lambda_{i_1}\cdots \lambda_{i_k} \left(\sum_{p=1}^{k} \frac{\tilde{U}_{i_p i_p}}{\lambda_{i_p}}  \right)\right]^2\\
		&+\frac{1}{\sigma_k(\lambda)}\sum_{i_1,\cdots, i_k}\lambda_{i_1}\cdots \lambda_{i_k}  \left[ \sum_{p=1}^{k} \frac{\tilde{U}_{i_p i_p}}{\lambda_{i_p}}  \right]^2\\
			&-\frac{1}{\sigma_k (\lambda)}  \sum_{i_1,\cdots, i_k}\lambda_{i_1}\cdots \lambda_{i_k}  \sum_{p=1}^{k} \left(\frac{\tilde{U}_{i_p i_p}}{\lambda_{i_p}}  \right)^2.
	\end{split}
\end{equation}
We show now that $d^{2} F_{k}(D,\widetilde{U})$ coincides indeed with the second-order differential $d^{2} f_{k}(\lambda,\xi)$ of the usual function $f_k(\lambda)=\log (\sigma_k(\lambda)),$ considered on $\Gamma_n$ by appropriating choosing $\lambda$ and $\xi.$

Indeed, we have for $\lambda\in \Gamma_{n}, \xi=(\xi_1, \cdots, \xi_n)\in \mathbb{R}^n:$ 
\begin{equation}\label{3.51}
	\begin{split}
		&  	d^{2} f_{k}(\lambda,\xi)=-\frac{1}{\sigma_k^2(\lambda)}\left[ \sum_{i_1,\cdots, i_k} \lambda_{i_1}\cdots \lambda_{i_k} \left(\sum_{p=1}^{k} \frac{\xi_{i_p}}{\lambda_{i_p}}  \right)\right]^2\\
		&+\frac{1}{\sigma_k(\lambda)}\sum_{i_1,\cdots, i_k}\lambda_{i_1}\cdots \lambda_{i_k}  \left[ \sum_{p=1}^{k} \frac{\xi_{i_p}}{\lambda_{i_p}}  \right]^2-\frac{1}{\sigma_k (\lambda)}  \sum_{i_1,\cdots, i_k}\lambda_{i_1}\cdots \lambda_{i_k}  \sum_{p=1}^{k} \left(\frac{\xi_{i_p}}{\lambda_{i_p}}  \right)^2.
	\end{split}
\end{equation}
From \eqref{3.50}, \eqref{3.51} it follows that
\begin{proposition}\label{md3.6}
Suppose
\begin{equation}\label{3.52}
	f_k(\lambda)=\log (\sigma_k(\lambda)),
\end{equation}	
considered on the positive cone $\Gamma_{n}$ of  $ \mathbb{R}^{n}.$ Suppose 
$$D={\rm diag}\left(\lambda_{1}, \lambda_{2}, \cdots, \lambda_{n}\right)>0,\quad \widetilde{U}={\rm diag}\left(\widetilde{U}_{11}, \widetilde{U}_{22}, \cdots, \widetilde{U}_{nn}\right).$$
Then we have 
\begin{equation*}\label{3.53}
	d^{2} F_{k}(D, \widetilde{U})=d^{2} f_{k}(\lambda, \xi)
\end{equation*}
if the following conditions fulfill
\begin{equation*}\label{3.54}
	\lambda=(\lambda_1, \cdots, \lambda_n)=\lambda(D),\quad 
	\xi=(\xi_1, \cdots, \xi_n)=(\widetilde{U}_{11}, \cdots,  \widetilde{U}_{nn}).
\end{equation*}
\end{proposition}
In next section we show that for  $2\le k\le n-1$ the function $f_{k}(\lambda)$ is strictly concave on a subset of $\Gamma_{n},$ i.e. its second-order differential is negative definite quadratic form on $\mathbb{R}^n_{\xi}$ when $\lambda$ varies on that subset.
\section{The strict concavity of the function $f_k(\lambda)$}
We consider the function $f_k(\lambda)=\log (\sigma_k(\lambda))$ on $\Gamma_{n}.$ For the index $i_1i_2\cdots i_k$ with $1\le i_1 < i_2< \cdots < i_k\le n$ we denote 
\begin{equation}\label{4.1}
	\gamma_{i_1\cdots i_k}=\frac{ \lambda_{i_1} \lambda_{i_2}\cdots \lambda_{i_k}}{\sigma_{k}(\lambda)}.
\end{equation}
It is obvious that 
\begin{equation}\label{4.2}
	0< \gamma_{i_1\cdots i_k}<1,
\end{equation}
\begin{equation}\label{4.3}
\sum_{i_1,\cdots, i_k}\gamma_{i_1\cdots i_k}=1.
\end{equation}
From \eqref{3.51} and \eqref{4.1} we rewrite $d^{2} f_{k}(\lambda, \xi)$ as follows 
\begin{equation}\label{4.4}
		\begin{split}
		  d^{2} f_{k}(\lambda, \xi)	=&-\left[ \sum_{i_1,\cdots, i_k} \gamma_{i_1\cdots i_k} \left(\sum_{p=1}^{k} \frac{\xi_{i_p}}{\lambda_{i_p}}  \right)\right]^2\\
		& +\sum_{i_1,\cdots, i_k} \gamma_{i_1\cdots i_k} \left[\sum_{p=1}^{k} \frac{\xi_{i_p}}{\lambda_{i_p}}  \right]^2-\sum_{i_1,\cdots, i_k} \gamma_{i_1\cdots i_k} \sum_{p=1}^{k} \left(\frac{\xi_{i_p}}{\lambda_{i_p}}  \right)^2.
	\end{split}
\end{equation}
Suppose $\eta=(\eta_1, \eta_2,\cdots, \eta_n)\in \mathbb{R}^n.$ For $\lambda=(\lambda_1, \cdots, \lambda_n)\in \Gamma_{n}$ we denote
\begin{equation}\label{4.5}
	\frac{1}{\lambda}=\left(\frac{1}{\lambda_1}, \frac{1}{\lambda_2}, \cdots, \frac{1}{\lambda_n}\right),
\end{equation}
\begin{equation}\label{4.6}
	\lambda\eta=(\lambda_1\eta_1,\lambda_2\eta_2,\cdots,\lambda_n\eta_n),
\end{equation}
\begin{equation}\label{4.7}
	\frac{1}{\lambda} \eta =\left( \frac{\eta_1}{\lambda_{1}}, \frac{\eta_2}{\lambda_{2}},\cdots, \frac{\eta_n}{\lambda_{n}} \right).
\end{equation}
It is obvious that $\frac{1}{\lambda}\in \Gamma_{n}$ and $\lambda\eta, 	\frac{1}{\lambda}\eta$ are linear transforms of the variable $\eta$ in $\mathbb{R}^n.$  

 We consider the quadratic form $  \widetilde{ d^{2} f_{k}(\lambda, \eta)}$ defined as follows
 \begin{equation}\label{4.8}
 	 \widetilde{ d^{2} f_{k}(\lambda, \eta)}=d^{2} f_{k}(\lambda, \lambda\eta).
 \end{equation}
From \eqref{4.4} and \eqref{4.8} we have 
\begin{equation}\label{4.9}
	\begin{split}
	 \widetilde{ d^{2} f_{k}(\lambda, \eta)}=	& -\left[ \sum_{i_1,\cdots, i_k} \gamma_{i_1\cdots i_k} \left(\sum_{p=1}^{k} \eta_{i_p}  \right)\right]^2\\
		& +\sum_{i_1,\cdots, i_k}\gamma_{i_1\cdots i_k} \left[\sum_{p=1}^{k} \eta_{i_p}  \right]^2-\sum_{i_1,\cdots, i_k}\gamma_{i_1\cdots i_k}\left[\sum_{p=1}^{k} \eta_{i_p}^2  \right].
	\end{split}
\end{equation}
\subsection{Simplifying  $\widetilde{ d^{2} f_{k}(\lambda, \eta)}$}
\begin{proposition}\label{md4.1}
		The quadratic form  $\widetilde{ d^{2} f_{k}(\lambda, \eta)}$ can be simplified as follows:
	\begin{equation}\label{4.10}
		  \widetilde{ d^{2} f_{k}(\lambda, \eta)}=-\frac{1}{\sigma_k^2(\lambda)}   \left[ \sum_{i=1}^n (\lambda_{i}\sigma_{k-1}^{(i)} (\lambda))^2 \eta^2_i+\sum_{i\ne j}\lambda_{i}\lambda_{j}a^{(k)}_{ij}(\lambda)\eta_i\eta_j \right]	,				
	\end{equation}
where 
\begin{equation}\label{4.11}
\begin{split}
		a^{(k)}_{ij}(\lambda) &=\frac{1}{2}\left[\sigma_{k-1}^{(i,j)} (\lambda)   \right]^2-\sigma_{k}^{(i,j)} (\lambda)  \sigma_{k-2}^{(i,j)} (\lambda),\\
	\sigma_{k-1}^{(i)} (\lambda)&=\left.\sigma_{k-1}(\lambda)\right|_{\lambda_i=0},\\
	\sigma_{k-1}^{(i,j)} (\lambda)&=\left.\sigma_{k-1}(\lambda)\right|_{\lambda_i=\lambda_j=0}.
\end{split}
\end{equation}	
	\end{proposition}
\begin{proof}
We consider the sum of the first and second  terms on the right-hand side of \eqref{4.9}.  From \eqref{4.3} we can rewrite this sum as follows
\begin{equation}\label{4.12}
		\begin{split}
		 -\left[ \sum_{i_1,\cdots, i_k} \gamma_{i_1\cdots i_k}   \left(\sum_{p=1}^{k} \eta_{i_p}  \right)\right]^2&+ \sum_{i_1,\cdots, i_k} \gamma_{i_1\cdots i_k} \left[\sum_{p=1}^{k}\eta_{i_p}\right]^2\\
		 &=\frac{1}{2}\sum_{i_1,\cdots, i_k\atop j_1,\cdots, j_k} \gamma_{i_1\cdots i_k}  \gamma_{j_1\cdots j_k} \left( \sum_{p=1}^{k} \eta_{i_p}-  \sum_{q=1}^{k} \eta_{j_q}\right)^2.
		\end{split}
	\end{equation}
We note that 
$$\sum_{p=1}^{k} \eta_{i_p}=\sum_{m=1}^{n} \eta_{m}-\sum_{\ell=1\atop i'_{\ell}\notin \{i\}_k}^{n-k} \eta_{i'_{\ell}},\quad 
\sum_{q=1}^{k} \eta_{j_q}=\sum_{m=1}^{n} \eta_{m}-\sum_{h=1\atop j'_h\notin \{j\}_k}^{n-k} \eta_{j'_h},$$
where $\{i\}_k$ stands for $\{i_1,\cdots, i_k\},$ $\{j\}_k$ for $\{j_1,\cdots, j_k\},$
$$\sum_{\ell=1\atop i'_{\ell}\notin \{i\}_k}^{n-k} \eta_{i'_{\ell}}= \sum_{  i'_{\ell}\in \{j\}_k\backslash\{i\}_k} \eta_{i'_{\ell}} +\sum_{i'_{\ell}\notin \{i\}_k\cup\{j\}_k}\eta_{i'_{\ell}}.$$
So we have
\begin{equation}\label{4.13}
	\sum_{p=1 }^{k} \eta_{i_p}- 	\sum_{q=1 }^{k} \eta_{j_q}=\sum_{  j'_h\in \{i\}_k\backslash\{j\}_k} \eta_{j'_k}-\sum_{  i'_{\ell}\in \{j\}_k\backslash\{i\}_k} \eta_{i'_{\ell}}.
\end{equation}
From \eqref{4.13} it follows that
\begin{equation*}\label{4.14}
	\begin{split}
	 \left(\sum_{p=1 }^{k} \eta_{i_p}- 	\sum_{q=1 }^{k} \eta_{j_q} \right)^2=&\left[ \sum_{  j'_h\in \{i\}_k\backslash\{j\}_k} \eta_{j'_h}\right]^2+\left[\sum_{  i'_{\ell}\in \{j\}_k\backslash\{i\}_k} \eta_{i'_{\ell}} \right]^2\\
	 &\hspace{1cm}- \left[ \sum_{  j'_h\in \{i\}_k\backslash\{j\}_k} \eta_{j'_h}\right]  \left[\sum_{  i'_{\ell}\in \{j\}_k\backslash\{i\}_k} \eta_{i'_{\ell}}  \right]\\
	 = &\sum_{  j'_h\in \{i\}_k\backslash\{j\}_k} (\eta_{j'_h})^2 +
	 \sum_{  j'_{\ell}, j'_h\in \{i\}_k\backslash\{j\}_k} \eta_{j'_{\ell}}\eta_{j'_h} +\sum_{  i'_{\ell}\in \{j\}_k\backslash\{i\}_k} (\eta_{i'_{\ell}})^2\\
	 &   +\sum_{i'_h, i'_{\ell}\in \{j\}_k\backslash\{i\}_k} \eta_{i'_h}\eta_{i'_{\ell}} -   \sum_{    j'_h\in \{i\}_k\backslash\{j\}_k\atop  i'_{\ell}\in \{j\}_k\backslash\{i\}_k} (\eta_{j'_h}\eta_{i'_{\ell}})
	\end{split}
\end{equation*}
Therefore 
\begin{equation*}\label{4.15}
	\begin{split}
	 \frac{1}{2}& \sum_{i_1,\cdots, i_k\atop j_1,\cdots, j_k } \gamma_{i_1\cdots i_k} \gamma_{j_1\cdots j_k} \left(\sum_{p=1}^k \eta_{i_p}- \sum_{q=1}^k \eta_{j_q} \right)^2\\
	&=\sum_{  i_1,\cdots, i_k\atop j_1,\cdots, j_k } \gamma_{i_1\cdots i_k} \gamma_{j_1\cdots j_k} \left[\sum_{  j'_h\in \{i\}_k\backslash\{j\}_k}  (\eta_{j'_{\ell}})^2  \right]+ \sum_{  i_1,\cdots, i_k\atop j_1,\cdots, j_k } \gamma_{i_1\cdots i_k} \gamma_{j_1\cdots j_k} \left[\sum_{  j'_{\ell}, j'_h\in \{i\}_k\backslash\{j\}_k}  \eta_{j'_{\ell}} \eta_{j'_h}  \right] \\
	&\hspace{3cm} -\frac{1}{2} \sum_{  i_1,\cdots, i_k\atop j_1,\cdots, j_k } \gamma_{i_1\cdots i_k} \gamma_{j_1\cdots j_k} \left[\sum_{ j'_h\in \{i\}_k\backslash\{j\}_k\atop i'_{\ell}\in \{j\}_k\backslash\{i\}_k}  \eta_{j'_h} \eta_{i'_{\ell}}  \right]\\
	&= \sum_{\ell=1}^{n}\left[ \sum_{\ell\in \{i\}_k} \left( \sum_{\ell\notin \{j\}_k}\gamma_{j_1\cdots j_k}  \right)\gamma_{i_1\cdots i_k} \right]\eta^2_{\ell}\\
	&\hspace{1cm} +\sum_{\ell,h=1\atop \ell\ne h}^{n} \left[-\frac{1}{2} \left(\sum_{i_1,\cdots,i_k \atop h\in \{i\}_k,   \ell\notin\{i\}_k} \gamma_{i_1\cdots i_k} \right) \left( \sum_{j_1,\cdots, j_k \atop \ell\in  \{j\}_k, h\notin \{j\}_k} \gamma_{j_1\cdots j_k}  \right)\right.\\
&\hspace{3cm} \left. +  \left(\sum_{ i_1,\cdots,i_k\atop \ell,h\in \{i\}_k} \gamma_{i_1\cdots i_k} \right) \left( \sum_{j_1,\cdots,j_k\atop   \ell,h\notin  \{j\}_k } \gamma_{j_1\cdots j_k}  \right)  \right]\eta_{\ell}\eta_h.
	\end{split}
\end{equation*}
We set 
\begin{equation*}\label{4.16}
	\sum_{j_1,\cdots,j_k\atop \ell\notin \{j\}_k} \gamma_{j_1\cdots j_k}=\gamma_{\ell}.
\end{equation*}
Then we have 
\begin{equation*}\label{4.17}
	\sum_{i_1,\cdots,i_k\atop \ell\in \{i\}_k} \left(\sum_{j_1,\cdots,j_k \atop \ell\notin \{j\}_k}\gamma_{j_1\cdots j_k} \right)\gamma_{i_1\cdots i_k}=\gamma_{\ell} 	\sum_{i_1,\cdots,i_k\atop \ell\in \{i\}_k}\gamma_{i_1\cdots i_k}.
\end{equation*}
Since $\gamma_{i_1\cdots i_k}=\frac{\lambda_{i_1} \cdots \lambda_{i_k}}{\sigma_{k}(\lambda)},$ then if $\ell\ne h$ we have
\begin{equation}\label{4.18}
 \left(\sum_{i_1,\cdots,i_k \atop h\in \{i\}_k,  \ell\notin \{i\}_k } \gamma_{i_1\cdots i_k} \right) \left(\sum_{j_1,\cdots,j_k\atop \ell\in \{j\}_k , h\notin  \{j\}_k } \gamma_{j_1\cdots j_k} \right)= \frac{\lambda_h \sigma_{k-1}^{(h,\ell)}(\lambda) \lambda_{\ell} \sigma_{k-1}^{(\ell,h)}(\lambda)}{\sigma_{k}^{2}(\lambda)},
\end{equation}
\begin{equation}\label{4.19}
	\left(\sum_{i_1,\cdots,i_k\atop \ell,h\in \{i\}_k } \gamma_{i_1\cdots i_k} \right) \left(\sum_{j_1,\cdots,j_k\atop \ell,h\notin \{j\}_k } \gamma_{j_1\cdots j_k} \right)= \frac{\lambda_{\ell}\lambda_h \sigma_{k-2}^{(\ell,h)}(\lambda) \sigma_{k}^{(\ell,h)}(\lambda)}{\sigma_{k}^{2}(\lambda)}.
\end{equation}
On the other hand 
\begin{equation*}\label{4.20}
\sum_{i_1,\cdots, i_k} \gamma_{i_1\cdots i_k}	\left[\sum_{p=1}^k \eta^2_{i_p} \right]=\sum_{\ell=1}^{n} \left(\sum_{i_1,\cdots,i_k,\atop \ell\in \{i\}_k }\gamma_{i_1\cdots i_k}	 \right) \eta^2_{\ell}=\sum_{\ell=1}^{n} (1-\gamma_{\ell}) \eta^2_{\ell},
\end{equation*}
and we have also 
\begin{equation}\label{4.21}
\begin{split}
	\gamma_{\ell}(1-	\gamma_{\ell})-(1-	\gamma_{\ell})&=-(1-	\gamma_{\ell})^2=- \left[ \sum_{\ell\in \{i\}_k }\gamma_{i_1\cdots i_k}	\right]^2\\
	&=-\frac{1}{\sigma_{k}^{2}(\lambda)} \left[\lambda_{\ell} \sigma_{k-1}^{(\ell)}(\lambda) \right]^2.
\end{split}
\end{equation}
From \eqref{4.9}, \eqref{4.12}, \eqref{4.18}-\eqref{4.21},  \eqref{4.10} follows.  
	\end{proof}
\subsection{The strict negative definiteness of  $\widetilde{d^{2} f_{k}(\lambda, \eta)}$}
From \eqref{4.8} we have the following 
\begin{remark}
	Suppose $\lambda\in \Gamma_{n}.$ Then the following assertions are true:
\begin{itemize}
	\item[(i)] $\widetilde{ d^{2} f_{k}(\lambda, \eta)}\le 0$ for all $\eta\in \mathbb{R}^n$ if and only if $ d^{2} f_{k}(\lambda, \xi)\le 0$ for all $\xi\in \mathbb{R}^n;$ 
	\item[(ii)] $\widetilde{ d^{2} f_{k}(\lambda, \eta)} < 0$ for all $\eta\in \mathbb{R}^n, \eta\ne 0$ if and only if $ d^{2} f_{k}(\lambda, \xi)\le 0$ for all $\xi\in \mathbb{R}^n, \xi\ne 0.$ 
	\end{itemize}
	\end{remark}
\begin{proposition}\label{bd4.1}
 The following assertion is true 	for all $k:$ $2\le k\le n-1,$  
 \begin{equation*}\label{4.23}
 	\widetilde{ d^{2} f_{k}(\lambda, \eta)} \le 0 \text{  for all } \lambda\in \Gamma_n, \eta \in \mathbb{R}^n.
 \end{equation*}
 \end{proposition}
\begin{proof}
We denote the first, second and third terms on the right-hand side of \eqref{4.4} by $(-A(\lambda, \xi)),$ $(B(\lambda, \xi))$   and $(-C(\lambda, \xi))  $ respectively. So 
\begin{equation}\label{4.24}
	d^{2} f_{k}(\lambda, \xi)=-A(\lambda, \xi)+  B(\lambda, \xi)- C(\lambda, \xi).
\end{equation}
We consider the function 
\begin{equation*}\label{4.25}
	g_k(\lambda)=\sqrt[k]{\sigma_{k}(\lambda)}
\end{equation*}
on the positive cone $\Gamma_{n} \subset \mathbb{R}^{n}.$ By calculating we have
\begin{equation*}\label{4.26}
		d^{2} g_{k}(\lambda, \xi)=\frac{(k-1)}{k^2}\left[\sigma_{k}(\lambda) \right]^{\frac{1}{k}} \left[-A(\lambda, \xi)+  B(\lambda, \xi)- C(\lambda, \xi)+\frac{1}{(k-1)} E(\lambda, \xi)\right],
\end{equation*}
where 
\begin{equation*}\label{4.27}
	E(\lambda, \xi)=\sum_{i_1,\cdots, i_k }\gamma_{i_1\cdots i_k}	\sum_{p,q=1\atop p\ne q}^k\frac{\xi_{i_p}\xi_{i_q}}{\lambda_{i_p}\lambda_{i_q} }.
\end{equation*}
It is well-known (\cite{2}) that 
\begin{equation}\label{4.28}
		d^{2} g_{k}(\lambda, \xi)\le 0, \text{ for all } \lambda\in \Gamma_n, \xi \in \mathbb{R}^n.
\end{equation}
We denote for $\eta\in \mathbb{R}^n, \lambda\in \Gamma_n$
\begin{equation*}\label{4.29}
	\widetilde{ d^{2} g_{k}(\lambda, \eta)} =	d^{2} g_{k}(\lambda, \lambda\eta).
\end{equation*}
From \eqref{4.28} it follows that 
\begin{equation*}\label{4.30}
		\widetilde{ d^{2} g_{k}(\lambda, \eta)}\le 0 \text{ for all } \lambda\in \Gamma_n, \eta \in \mathbb{R}^n.
\end{equation*}
We have
\begin{equation}\label{4.31}
\widetilde{ d^{2} g_{k}(\lambda, \eta)}	=\frac{(k-1)}{k^2}\left[\sigma_{k}(\lambda) \right]^{\frac{1}{k}} \left[-\widetilde{A(\lambda, \eta)}+  \widetilde{B(\lambda, \eta)}- \widetilde{C(\lambda, \eta)}+\frac{1}{(k-1)} \widetilde{E(\lambda, \eta)}\right],
\end{equation}
where 
\begin{equation*}\label{4.32}
	\widetilde{A(\lambda, \eta)}=\left[ \sum_{i_1,\cdots, i_k} \gamma_{i_1\cdots i_k}	\left(\sum_{p=1}^k \eta_{i_p} \right) \right]^2,
\end{equation*}
\begin{equation*}\label{4.33}
\widetilde{B(\lambda, \eta)}= \sum_{i_1,\cdots, i_k} \gamma_{i_1\cdots i_k}	\left[\sum_{p=1}^k \eta_{i_p}   \right]^2,
\end{equation*}
\begin{equation*}\label{4.34}
	 \widetilde{C(\lambda, \eta)} = \sum_{i_1,\cdots, i_k} \gamma_{i_1\cdots i_k}	 \sum_{p=1}^k \eta_{i_p}^2,
\end{equation*}
\begin{equation}\label{4.35}
	\widetilde{E(\lambda, \eta)} = \sum_{i_1,\cdots, i_k} \gamma_{i_1\cdots i_k}	  
	\sum_{p,q=1\atop p\ne q}^k \eta_{i_p}\eta_{i_q}.
\end{equation}
We note that $(-\widetilde{A(\lambda, \eta)})\le 0$ and 
\begin{equation}\label{4.36}
	\widetilde{B(\lambda, \eta)}-\widetilde{C(\lambda, \eta)}=\widetilde{E(\lambda, \eta)}.
\end{equation}
Since  $\widetilde{ d^{2} g_{k}(\lambda, \eta)}\le 0,$ from \eqref{4.31},  \eqref{4.36} it follows that 
$$-\widetilde{A} +\frac{k}{(k-1)}\widetilde{E} \le 0$$
and 
\begin{equation*}\label{4.37}
	\widetilde{E}\le \frac{(k-1)}{k}\widetilde{A}.
\end{equation*}
From \eqref{4.24} we have for any $\lambda\in \Gamma_n$ and $\eta\in \mathbb{R}^n,$
\begin{equation}\label{4.38}
	\begin{split}
		\widetilde{ d^{2} f_{k}(\lambda, \eta)}&=- 	\widetilde{A} +\widetilde{B} -\widetilde{C} =-\widetilde{A} +\widetilde{E} \\
		&\le -\widetilde{A} +\frac{k-1}{k}\widetilde{A} =-\frac{1}{k}\widetilde{A} \le 0.
	\end{split}
\end{equation}
	\end{proof}
\begin{proposition} \label{md4.2}
	Suppose $2\le k\le n-1$ and $\lambda\in \Gamma_{n}.$ Then the following assertion is true:
	\begin{equation*}\label{4.39}
	\widetilde{ d^{2} f_{k}(\lambda, \eta)}<0 \text{ for all } \eta\in \mathbb{R}^n, \eta\ne 0		
	\end{equation*}
if and only if the matrix
	\begin{equation}\label{4.40}
		G_k(\lambda)=\begin{bmatrix}
			0&\sigma_{k-2}^{(1,2)} (\lambda)&\cdots&\sigma_{k-2}^{(1,n)}(\lambda)\\ 
			\sigma_{k-2}^{(1,2)}(\lambda)&0&\cdots&\sigma_{k-2}^{(2,n)}(\lambda) \\ 
			\cdots&\cdots&\cdots&\cdots\\ 
			\sigma_{k-2}^{(1,n)}(\lambda) &\sigma_{k-2}^{(2,n)}(\lambda)&\cdots&0\\ 
		\end{bmatrix}_{n\times n}
	\end{equation}
is not degenerative.
\end{proposition}
\begin{proof}
	Suppose $\lambda\in \Gamma_n$ and there exists $\overline{\eta}\in \mathbb{R}^n$ such that
	\begin{equation}\label{4.41}
	\widetilde{ d^{2} f_{k}(\lambda, \overline{\eta})}=0.	
	\end{equation}
From \eqref{4.41} and \eqref{4.38}   it follows that 	
\begin{equation}\label{4.42}
\widetilde{A(\lambda,\overline{\eta})} 	=0,
\end{equation}
\begin{equation}\label{4.43}
	\widetilde{E(\lambda,\overline{\eta})} 	=0.
\end{equation}
Since 
\begin{equation}\label{4.43'}
-\widetilde{A(\lambda,\eta)} + 	\widetilde{E(\lambda,\eta)}\le 0, \text{ for all } \eta\in \mathbb{R}^n, 	
\end{equation}
\begin{equation}\label{4.27'}
\widetilde{A(\lambda,\eta)} \ge 0, \text{ for any } \eta\in \mathbb{R}^n,	
\end{equation}
then from \eqref{4.42}-\eqref{4.27'} it follows that 
\begin{equation}\label{4.44}
	-\frac{\partial}{\partial \eta_j} \widetilde{A(\lambda,\overline{\eta})}+ \frac{\partial}{\partial \eta_j} \widetilde{E(\lambda,\overline{\eta})}=0, \quad j=1,\cdots, n,
\end{equation}
\begin{equation}\label{4.45}
	-\frac{\partial}{\partial \eta_j} \widetilde{A(\lambda,\overline{\eta})}=0, \quad j=1,\cdots, n.
\end{equation}
From \eqref{4.44},  \eqref{4.45}  we have 
\begin{equation}\label{4.46}
	\frac{\partial}{\partial \eta_j} \widetilde{E(\lambda,\overline{\eta})}=0, \quad j=1,\cdots, n.
\end{equation}
From \eqref{4.35}  it follows that 
\begin{equation*}\label{4.47}
	\widetilde{E(\lambda,\eta)}=\sum_{i=1\atop i\ne j}^n \left(\sum_{i,j \in \{i\}_k} \gamma_{i_1\cdots i_k} \right)\eta_{i}\eta_{j}=\frac{1}{\sigma_k(\lambda)}\sum_{i,j =1\atop i\ne j}^n \lambda_i \lambda_j\sigma_{k-2}^{(i,j)}(\lambda) \eta_i \eta_j.
\end{equation*}
\begin{equation}\label{4.48}
	\frac{\partial}{\partial \eta_j} \widetilde{E(\lambda,\eta)} =\frac{1}{\sigma_k(\lambda)}\sum_{i =1\atop i\ne j}^n \lambda_i \lambda_j\sigma_{k-2}^{(i,j)}(\lambda) \eta_i.
\end{equation}
From \eqref{4.46}, \eqref{4.48}  it follows that $\overline{\eta}=0$ if and only if the matrix $G_{k}(\lambda),$ defined by \eqref{4.40}, is not degenerative.
	\end{proof}
\begin{proposition}\label{md4.3}
  Suppose $ 3 \leq k \leq n-1,$ $\lambda\in \Gamma_n$ and $\widetilde{ d^{2} f_{k}(\lambda, \eta)}<0,$ for  all $\eta\in \mathbb{R}^n, \eta\ne 0.$ Then the following assertion is true:
  \begin{equation}\label{4.52}
  	\widetilde{ d^{2} f_{n-k+2}(\frac{1}{\lambda}, \eta)}<0  \text{ for  all } \eta\in \mathbb{R}^n, \eta\ne 0.
  \end{equation}
\end{proposition}
\begin{proof}
	From the equality
	\begin{equation*}\label{4.53}
		\sigma_{k}(\lambda)=(\lambda_{1}\lambda_{2}\cdots \lambda_{n}) 	\sigma_{n-k}\left(\frac{1}{\lambda}\right), \lambda\in \Gamma_n,
	\end{equation*}
it follows that for all $i\ne j$ and $\lambda\in \Gamma_n$
\begin{equation*}\label{4.54}
	\sigma_{k-2}^{(i,j)}(\lambda)=\frac{\lambda_{1}\lambda_{2}\cdots \lambda_{n}}{\lambda_{i}\lambda_j}\sigma_{(n-2)-(k-2)}^{(i,j)}\left(\frac{1}{\lambda}\right)= \frac{\lambda_{1}\lambda_{2}\cdots \lambda_{n}}{\lambda_{i}\lambda_j}\sigma_{(n-k+2)-2}^{(i,j)}\left(\frac{1}{\lambda}\right),
\end{equation*}
from which we obtain 
\begin{equation*}\label{4.55}
\sigma_{(n-k+2)-2}^{(i,j)}\left(\frac{1}{\lambda}\right)=	\frac{\lambda_{i}\lambda_j}{\lambda_{1}\lambda_{2}\cdots \lambda_{n}}\sigma_{n-2}^{(i,j)}(\lambda).
\end{equation*}
Hence we have
$$\det\left[ G_{n-k+2}\left(\frac{1}{\lambda}\right)\right]=\frac{1}{(\lambda_{1}\lambda_{2}\cdots \lambda_{n})^{n-2}}\det[G_k(\lambda)],$$
from which it follows that if $\det[G_k(\lambda)]\ne 0$ then 
$$\det\left[G_{n-k+2}(\frac{1}{\lambda})\right] \ne 0,$$
  and \eqref{4.52} follows from Proposition \ref{md4.2}.
	\end{proof}
   \begin{proposition}\label{md4.4}
  	Suppose $k\in \{ 2,3,4,n-2,n-1\}$ and $\lambda\in \Gamma_{n}.$ Then the following assertion is true:
  	$$\widetilde{ d^{2} f_{k}(\lambda, \eta)}<0\text{ for all } \eta\in \mathbb{R}^n, \eta\ne 0.$$
  \end{proposition}
  \begin{proof}
  	(i)  When $ k=n-1,$   then  $ \sigma_{k}^{(i, j)}(\lambda)=0$ if $i\ne j$  and 
  	$$a_{ij}^{(n-1)}(\lambda) =\frac{1}{2}\left[ \sigma_{n-2}^{(i, j)}(\lambda)\right]^2.$$
  	From  \eqref{4.10} and the fact that $\sigma_{n-2}^{(i)}(\lambda)=\sum_{j\ne i}\sigma_{n-2}^{(i,j)}(\lambda)$ it follows that
  	\begin{equation}\label{4.50}
  		\begin{split}
  			\widetilde{ d^{2} f_{n-1}(\lambda, \eta)}&=-\frac{1}{ \sigma_{n-1}^{2}(\lambda)
  			}\left[\sum_{i=1}^{n} (\lambda_i \sigma_{n-2}^{(i)}(\lambda))^2 \eta^2_i    +\frac{1}{2}\sum_{i\ne j} \lambda_{i} \lambda_{j} \left[ \sigma_{n-2}^{(i, j)}(\lambda) \right]^2 \eta_i\eta_j 
  			\right]	\\
  			& \le -\frac{1}{ \sigma_{n-1}^{2}(\lambda)
  			} \left[\frac{1}{2}\sum_{i=1}^{n} (\lambda_{i} \sigma_{n-2}^{(i)}(\lambda))^2\eta_i^2+ \frac{1}{2}\sum_{ i\ne j}^n\lambda_{i}^2 \left[ \sigma_{n-2}^{(i, j)}(\lambda)\right]^2\eta_i^2\right.\\
  			&\hspace{5cm}\left.+\frac{1}{2}\sum_{i\ne j} \lambda_i\lambda_j \left[\sigma_{n-2}^{(i, j)}(\lambda) \right]^2\eta_i\eta_j \right]\\
  			& = -\frac{1}{2 \sigma_{n-1}^{2}(\lambda)
  			}  \left[\sum_{i=1}^n \left[\lambda_{i} \sigma_{n-2}^{(i)}(\lambda)\right]^2\eta_i^2 +\frac{1}{2} \sum_{ i\ne j}(\lambda_i\eta_{i}+\lambda_j\eta_{j} )^2\left[\sigma_{n-2}^{(i, j)}(\lambda) \right]^2\right].
  		\end{split}	
  	\end{equation}
  	So, it follows that $\widetilde{ d^{2} f_{n-1}(\lambda, \eta)}<0$ for all $\lambda\in \Gamma_n, $ $\eta\in \mathbb{R}^n, \eta\ne 0.$
  	
  	(ii)  When  $k=n-2,$ from \eqref{4.11}   we have 
  	\begin{equation*}\label{4.51}
  		a_{ij}^{(n-2)}(\lambda) =\frac{1}{2}  \sigma_{n-3}^{(i, j)}(\lambda^2),
  	\end{equation*}
  	where $\lambda^2:= (\lambda_{1}^2,\lambda_{2}^2,\cdots, \lambda_{n}^2).$
  	
  	From \eqref{4.10},  \eqref{4.50} and the fact that $\sigma_{n-3}^{(i)}(\lambda^2)=\frac{1}{2}\sum_{j\ne i}\sigma_{n-3}^{(i,j)}(\lambda^2)$ we have
  	\begin{equation*}\label{4.52T}
  		\begin{split}
  			&\widetilde{ d^{2} f_{n-1}(\lambda, \eta)}=-\frac{1}{ \sigma_{n-2}^{2}(\lambda)
  			}\left[\sum_{i=1}^{n} (\lambda_i \sigma_{n-3}^{(i)}(\lambda))^2 \eta^2_i    +  \frac{1}{2} \sum_{i\ne j} \lambda_{i} \lambda_{j} \sigma_{n-3}^{(i, j)}(\lambda^2)   \eta_i\eta_j 
  			\right]\\
  			&\le  -\frac{1}{ \sigma_{n-2}^{2}(\lambda)
  			} \left[ \sum_{i=1}^n \lambda_{i}^2 \Big[( \sigma_{n-3}^{(i)}(\lambda^2) +2\sum_{1\le i_1<\cdots < i_{n-4}\le n} \lambda_{i_1}^2\cdots \lambda_{i_{n-4}}^2 \sigma_2^{(i,\lambda_{i_1},\cdots,  \lambda_{i_{n-4}})} (\lambda) )\Big]\eta_i^2\right.\\
  			&\hspace{7cm}+\left. \frac{1}{2} \sum_{i\ne j} \lambda_{i}\lambda_j \sigma_{n-3}^{(i, j)}(\lambda^2)   \eta_i\eta_j\right]\\
  			&=-\frac{1}{ \sigma_{n-2}^{2}(\lambda)} \left[\frac{1}{4}\sum_{i\ne j}(\lambda_{i}^2\eta_{i}^2+\lambda_{j}^2\eta_{j}^2)\sigma_{n-3}^{(i,j)}(\lambda^2) \right.\\
  			&\quad +\frac{1}{2}\sum_{i\ne j}\lambda_i \lambda_j \sigma_{n-3}^{(i, j)}(\lambda^2)   \eta_i\eta_j
  			\left.+2 \sum_{i=1}^n\lambda_{i}^2 \Big[ \sum_{1\le i_1<\cdots< i_{n-4}\le n} \lambda_{i_1}^2\cdots \lambda_{i_{n-4}}^2 \sigma_2^{(i,i_1,\cdots,  i_{n-4})} (\lambda)\Big]\eta_i^2\right]\\
  			&\le - \frac{2}{ \sigma_{n-2}^{2}(\lambda)}\sum_{i=1}^n\left[\lambda_{i}^2\Big[\sum_{1\le i_1<\cdots< i_{n-4}\le n} \lambda_{i_1}^2\cdots \lambda_{i_{n-4}}^2 \sigma_2^{(i,i_1,\cdots,  i_{n-4})} (\lambda) \Big]\eta^2_i\right],
  		\end{split}
  	\end{equation*}
  	from which it follows that $\widetilde{ d^{2} f_{n-2}(\lambda, \eta)}<0$ for all $\lambda\in \Gamma_n,$ for all $\eta\in \mathbb{R}^n, \eta\ne 0.$
  	
  	(iii) The case $k=2$
  	
  	In this case $\sigma_{k-2}^{(i, j)}(\lambda)=\sigma_{0}^{(i, j)}(\lambda)=1, i\ne j,$ and 
  	$$\det G_2(\lambda)= \begin{vmatrix}
  		0&1&1&\cdots &1\\ 
  		1&0&1&\cdots&1\\ 
  		1&1&0&\cdots&1\\ 
  		\cdots&\cdots&\cdots&\cdots&\\ 
  		1&1&1&\cdots&0\\ 
  	\end{vmatrix}
  	=(-1)^{n-1}(n-1)\ne 0.$$
  	From Proposition \ref{md4.2} it follows that $\widetilde{ d^{2} f_{2}(\lambda, \eta)}<0$ for  all $\lambda\in \Gamma_n,$ and  $\eta\in \mathbb{R}^n, \eta\ne 0.$
  	
  	(iv) The cases $k=3$ and $k=4$ follow from Proposition \ref{md4.2} and the verified cases $k=n-1$ and $k=n-2$ respectively.
  \end{proof}
Suppose $0<\gamma_k<1.$ We introduce now a subset $\Sigma_{(\gamma_k)}$ of $\Gamma_{n}.$
\begin{definition}\label{dn4.1}
	The set $\Sigma_{(\gamma_k)}$ consists of all $\lambda=(\lambda_{1},\lambda_{2},\cdots,\lambda_{n})\in \Gamma_{n}$ such that 
	\begin{equation}\label{4.56}
		\lambda_{\min}\ge \gamma_k \lambda_{\max}
	\end{equation}
where $\lambda_{\min}=\min_{1\le j\le n}(\lambda_{j}),$ $\lambda_{\max}=\max_{1\le j\le n}(\lambda_{j}) $ and $\gamma_k$ satisfies the following conditions:

(i) If $k\in \{ 2,3,4,n-2,n-1\},$ then $\gamma_k$ is a some positive number that is less than 1;

(ii) If $[\frac{n}{2}]+1\le k \le n-3,$ then
\begin{equation}\label{4.57}
	\gamma_k=\frac{n-k}{k};
\end{equation}

(iii) If $5\le k\le \left[\frac{n}{2}\right],$ then 
$$\gamma_{k}=\gamma_{n-k+2}=\frac{k-2}{n-(k-2)}.$$
\end{definition}
\begin{remark}\label{nx4.2}
	We note that $\Sigma_{(\gamma_{k})}$ is a convex cone. Moreever, the following assertions are true:
	
	(i) If $\lambda\in \Sigma_{(\gamma_k)},$ then $\frac{1}{\lambda}\in  \Sigma_{(\gamma_k)};$
	
	(ii) If $[\frac{n}{2}]+1\le k_1<k_2\le n-3,$ then $\gamma_{k_1}>\gamma_{k_2}$ and
	$$ \Sigma_{(\gamma_{k_2})} \supset \Sigma_{(\gamma_{k_1})};$$
	
	(iii) If $5\le k_1<k_2\le [ \frac{n}{2}],$ then $\gamma_{k_1}<\gamma_{k_2}$ and 
		$$ \Sigma_{(\gamma_{k_2})} \subset \Sigma_{(\gamma_{k_1})};$$
		
		(iv) If $5\le k_1\le k_2=n-k_1+2,$ then 	$ \Sigma_{(\gamma_{k_1})} = \Sigma_{(\gamma_{k_2})}.$
	\end{remark} 
\begin{proposition}\label{md4.5}
Suppose $\left[\frac{n}{2}\right]+1 \leq  k \leq  n-3	$ and 
$\lambda=\left(\lambda_{1}, \lambda_{2},\cdots, \lambda_{n}\right) \in  \Sigma_{(\gamma_{k})}.$  Then the following assertions hold\\
(i) If $ \lambda_{1} \geq  \lambda_{2} \geq \cdots \geq  \lambda_{n}>0,$ then \begin{equation}\label{4.58}
	\lambda_{n} \geq \frac{\sigma_{k}^{(n)}(\lambda)}{\sigma_{k-1}^{(n)}(\lambda)};
\end{equation}
(ii) If $ \lambda_{1} \geq  \lambda_{2} \geq \cdots \geq  \lambda_{n}>0,$ then \begin{equation}\label{4.59}
	\frac{\sigma_{k}^{(n)}(\lambda)}{\sigma_{k-1}^{(n)}(\lambda)}\ge \frac{\sigma_{k}^{(j)}(\lambda)}{\sigma_{k-1}^{(j)}(\lambda)}, \quad j=1,2,\cdots n-1;
\end{equation}
(iii) In general, we have
\begin{equation}\label{4.60}
	\lambda_j\ge \frac{\sigma_{k}^{(j)}(\lambda)}{\sigma_{k-1}^{(j)}(\lambda)}, \quad j=1,2,\cdots n.
\end{equation}
	\end{proposition}
\begin{proof}
First we prove that 
\begin{equation}\label{4.61}
	\sup_{\lambda_{1}\ge \cdots \ge \lambda_{n-1}>0}\frac{\sigma_{k}^{(n)}(\lambda)}{\sigma_{k-1}^{(n)}(\lambda)}=\gamma_{k} \lambda_{1}.
\end{equation}	
	We consider the function
	$$g(\lambda')=g(\lambda_{1},\cdots, \lambda_{n-1})=\frac{\sigma_{k}^{(n)}(\lambda)}{\sigma_{k-1}^{(n)}(\lambda)}.$$
This function is concave on the convex set $\lambda_{1}\ge \cdots \ge \lambda_{n-1}>0$ so we have 
$$D^{2} g\left(\lambda_{1},\cdots, \lambda_{n-1}\right) \leq 0$$
and 
$$\frac{\partial^{2} g\left(\lambda_{1}, \cdots, \lambda_{n-1}\right)}{\partial \lambda_j^{2}} \leq  0, \quad j=1,2,\cdots, n-1.$$
Therefore, the function $\frac{\partial g\left(\lambda_{1}, \cdots \lambda_{n-1}\right)}{\partial \lambda_{j}},$ considered as a function of $\lambda_j,$ does not increase as $\lambda_j$ increases. We have 
$$\frac{\partial g\left(\lambda'\right)}{\partial \lambda_{j}}=\frac{\sigma_{k-1}^{(n, j)}(\lambda) \sigma_{k-1}^{(n)}(\lambda)-\sigma_{k}^{(n)}(\lambda) \cdot \sigma_{k-2}^{(n, j)}(\lambda)}{\left[\sigma_{k-1}^{(n)}(\lambda)\right]^{2}}$$
and 
$$\left.\frac{\partial g\left(\lambda'\right)}{\partial \lambda_{j}}\right|_{\lambda_{1}=\cdots=\lambda_{n-1}}=\frac{(n-2) !(n-1) !}{(n-k) !(n-k-1) !(k-1) ! k !}>0.$$
It follows that 
\begin{equation}\label{4.62}
	\frac{\partial g\left(\lambda'\right)}{\partial \lambda_{j}}>0, \text{  for all } \lambda_{1}\ge \lambda_{2}\ge \cdots \ge \lambda_{n-1}>0, j=1,2,\cdots, n-1,
\end{equation}
and from \eqref{4.57} the function $g\left(\lambda'\right)$ attains its maximum at $\lambda'=(\lambda_{1},\lambda_{1},\cdots, \lambda_{1})$ with the value $\gamma_{k} \lambda_{1},$ from which it follows \eqref{4.61}. From \eqref{4.62} the inequality \eqref{4.59} follows. The inequalities \eqref{4.60} are easy consequence of \eqref{4.58} and \eqref{4.59}.
	\end{proof}
\begin{proposition}\label{md4.6}
Suppose $\left[\frac{n}{2}\right]+1 \leq  k \leq  n-3, \lambda \in \Sigma_{\left(\gamma_{k}\right)}$
and
\begin{equation}\label{4.63}
	\widetilde{d^{2} f_{n-k}\left(\frac{1}{\lambda}, \eta\right) }< 0, \text { for all } \eta \in \mathbb{R}^{n}, \eta \neq 0.
\end{equation}
Then 
\begin{equation}\label{4.64}
\widetilde{	d^{2} f_{k}(\lambda, \eta)} < 0, \quad \text { for all }   \eta \in \mathbb{R}^{n}, \eta \neq 0.
\end{equation}
\end{proposition}

\begin{proof}
From the following equalities:
\begin{align*}
	\sigma_{n-k}^{(j)}\left(\frac{1}{\lambda}\right)&=\frac{1}{\left(\lambda_{1} \cdots \lambda_{n}\right)} \sigma_{k}^{(j)}(\lambda),   \\
\sigma_{n-k}\left(\frac{1}{\lambda}\right)&=\frac{1}{\left(\lambda_{1} \cdots \lambda_{n}\right)} \sigma_{k}(\lambda),\\
\frac{1}{\lambda_{i} \lambda_{j}}a_{ij}^{(n-k)} \left(\frac{1}{\lambda}\right)&=\frac{1}{\left(\lambda_{1} \cdots \lambda_{n}\right)} \lambda_{i} \lambda_{j} a_{ij}^{(k)}(\lambda),
\end{align*}
where $a_{i j}^{(k)}(\lambda)$ is defined by \eqref{4.11}, it follows from \eqref{4.63} that
\begin{equation}\label{4.65}
\widetilde{d^2f_{n-k}\left(\frac{1}{\lambda},\eta\right)}=\frac{-1}{\sigma_{k}^2(\lambda)}\left[ \sum_{j=1}^n \left[\sigma_{k}^{(j)}(\lambda)\right]^2\eta_j^2+ \sum_{i\ne j}\lambda_{i} \lambda_{j}a_{ij}^{(k)}(\lambda)\eta_i\eta_j\right]<0.
\end{equation}
So, from \eqref{4.60},  \eqref{4.65}  we have  
\begin{align*}
	\widetilde{d^2f_{k}\left(\lambda,\eta\right)}&=\frac{-1}{\sigma_{k}^2(\lambda)}\left[ \sum_{j=1}^n \left[\lambda_{j}\sigma_{k-1}^{(j)}(\lambda)\right]^2\eta_j^2+ \sum_{i\ne j}\lambda_{i} \lambda_{j}a_{ij}^{(k)}(\lambda)\eta_i\eta_j\right]  \\
	& \le  	\widetilde{d^2f_{n-k}\left(\frac{1}{\lambda},\eta\right)}<0\quad  \text{  for all } \eta\in \mathbb{R}^n, \eta\ne 0.
\end{align*}
	\end{proof}
\begin{proposition}\label{md4.8}
Suppose $5\le k\le n-3$ and $\lambda\in \Sigma_{(\gamma_{k})}.$ Then we have 
\begin{equation}\label{4.66}
		\widetilde{d^2f_{k}\left(\lambda,\eta\right)}>0 \quad \text{  for all }  \eta\in \mathbb{R}^n, \eta\ne 0.
\end{equation}	
	\end{proposition}
\begin{proof}
We prove \eqref{4.66} by induction with respect to a group of 4 parameters of $k:$
$$\{ 2h+3, 2h+4, n-(2h+1), n-(2h+2)\}$$
with $h=1,2,\cdots, \left[\frac{n}{4}\right]-2.$

First step. When $h=1$ then $k$ is one of the values $\{5,6,n-3,n-4\}.$ We show that if $\lambda\in \Sigma_{(\gamma_{k})},$ then \eqref{4.66} holds. 
Indeed, applying Proposition \ref{md4.6} with $k=n-3$ and $k=n-4$ respectively,  from Proposition \ref{md4.4} it follows that \eqref{4.66} holds for the cases $k=n-3$ and $k=n-4$ respectively.

Suppose $\lambda\in \Sigma_{(\gamma_{5})}.$ Since $\gamma_{5}=\gamma_{n-3},$ then $\lambda\in \Sigma_{(\gamma_{n-3})}$ and therefore  \eqref{4.66} holds when $k=n-3.$ From Proposition \ref{md4.3} it follows that \eqref{4.66} holds for $k=5.$ The case $k=6$ is proved by the same way.

Second step. Suppose \eqref{4.66} holds  when $k$ is one of the following values:
$$\{2h_0+3, 2h_0+4, n-(2h_0+1), n-(2h_0+2) \}$$
with $1\le h_0\le \left[\frac{n}{4}\right]-4.$ We prove that if $k$ is one of the values: 
$$\{2h_0+5, 2h_0+6, n-(2h_0+3), n-(2h_0+4) \},$$
then \eqref{4.66} holds if $\lambda\in \Sigma_{(\gamma_{k})},$ where  $1\le h_0\le \left[\frac{n}{4}\right]-4.$
Indeed, first we consider the case $k=n-\left(2 h_{0}+3\right).$ Suppose $\lambda,\frac{1}{\lambda}\in \Sigma_{(\gamma_{n-(2h_0+3)})}.$ Since 
$$\gamma_{n-(2h_0+3)}=\gamma_{n-[n-(2h_0+3)]+2}=\gamma_{2h_0+5}>\gamma_{2h_0+3},$$
then by Remark \ref{nx4.2} we have
$$\lambda, \frac{1}{\lambda}\in \Sigma_{(\gamma_{2h_0+3})}.$$
It follows from the assumptions of the induction, where \eqref{4.66} holds for $k=2h_0+3,$ that 
\begin{equation}\label{4.68}
		\widetilde{d^2f_{2h_0+3}\left(\lambda,\eta\right)}>0, \quad  \widetilde{d^2f_{2h_0+3}\left(\frac{1}{\lambda},\eta\right)}>0, \quad \eta\ne 0.
\end{equation}
Then from Proposition \ref{md4.6} and from \eqref{4.68} we have that 
\begin{equation}\label{4.69}
	\widetilde{d^2f_{n-(2h_0+3)}\left(\frac{1}{\lambda},\eta\right)}>0, \quad  \widetilde{d^2f_{n-(2h_0+3)}\left(\lambda,\eta\right)}>0, \quad \eta\ne 0.
\end{equation}
The case $k=n-(2h_0+4)$ is proved by the same way.

Consider now the case $k=2h_0+5.$ Suppose $\lambda, \frac{1}{\lambda}\in \Sigma_{(\gamma_{2h_0+5})}.$ Since 
$$\gamma_{2h_0+5}=\gamma_{n-(2h_0+5)+2}=\gamma_{n-(2h_0+3)}$$
we have that 
\begin{equation}\label{4.70}
	\lambda, \frac{1}{\lambda}\in \Sigma_{(\gamma_{n-(2h_0+3)})}.
\end{equation}
By the result, that has been proved just above, it follows that \eqref{4.69} holds. By applying Proposition \ref{md4.3}, from \eqref{4.69} we conclude that 
\begin{equation*} 
	\widetilde{d^2f_{2h_0+5}\left(\lambda,\eta\right)}>0, \quad  \widetilde{d^2f_{2h_0+5}\left(\frac{1}{\lambda},\eta\right)}>0, \quad \eta\ne 0.
\end{equation*}
The case $k=2h_0+6$ is proved by the same way.

Third step. We consider all the remaining cases of $k,$ obtained when $h_0=\left[\frac{n}{4}\right]-2.$ We have seen that when $1\le h\le \left[\frac{n}{4}\right]-3,$ then \eqref{4.66} holds already for the following cases of $k:$ 
\begin{equation}\label{4.67}
	\begin{split}
		5& \le k=2h+3\le 2 \left[\frac{n}{4}\right]-3, \\
	6	& \le k=2h+4\le 2 \left[\frac{n}{4}\right]-2, \\ 
	k&=n-(2h+1)\ge n-2\left[\frac{n}{4}\right]+5,\\
		k&=n-(2h+2)\ge n-2\left[\frac{n}{4}\right]+4.
	\end{split}
\end{equation}
So, we have to consider the following cases of $k$ when $h_0=\left[\frac{n}{4}\right]-2:$
\begin{equation*} 
	\begin{split}
			k&=2h_0+3= 2\left[\frac{n}{4}\right]-1,\\
			k&=2h_0+4= 2\left[\frac{n}{4}\right], \\ 
		k&=n-(2h_0+1) = n-2\left[\frac{n}{4}\right]+3,\\
	k&=n-(2h_0+2) = n-2\left[\frac{n}{4}\right]+2
	\end{split}
\end{equation*}
and the following two remaining cases of $k:$ 
\begin{equation*} 
	\begin{split}
		k&=n-2\left[\frac{n}{4}\right]+1,\\
		k&=n-2\left[\frac{n}{4}\right].
	\end{split}
\end{equation*}
By considering each of the cases $n=4m,$ $n=4m+1,$ $n=4m+2,$ $n=4m+3$ and applying Propositions \ref{md4.3}, \ref{md4.6} and the verified cases \eqref{4.67}, we can easily prove that  \eqref{4.66} holds for all these values of $k$ by the same argument, described in the second step.
	\end{proof}
\subsection{The uniform negative definiteness of  $\widetilde{ d^{2} f_{k}(\lambda, \eta)}$}
\begin{theorem}\label{dl4.1}
	Suppose $2 \le k \le n-1.$  Then there exists $\gamma^{(k)}>0$ such that for all $\lambda \in \Sigma_{(\gamma_k)},$  $ \eta \in \mathbb{R}^{n}:$
	\begin{equation}\label{4.58T}
		\widetilde{ d^{2} f_{k}(\lambda, \eta)}	\le -\gamma^{(k)}|\eta|^2.
	\end{equation}
	This means that the negative definiteness of the form $\widetilde{ d^{2} f_{k}(\lambda, \eta)}$  is uniform with respect to  $\lambda \in \Sigma_{(\gamma_k)}.$ 
\end{theorem}
\begin{proof}
Since $\gamma_{i_1i_2\cdots i_k},$ defined by \eqref{4.1},  is homogenous of degree $0$ with respect to $\lambda\in \Gamma_n,$ from \eqref{4.9} it follows that $\widetilde{ d^{2} f_{k}(\lambda, \eta)}$ is also homogenous of degree  $0$ with respect to $\lambda\in \Sigma_{(\gamma_k)}.$ So we consider $\widetilde{ d^{2} f_{k}(\lambda, \eta)}$ as defined on the set 
$$T_{(\gamma_k)}=\left(\Sigma_{(\gamma_k)}\cap S_{\lambda}^{n-1} \right)\times S_{\eta}^{n-1},$$
where $S_{\lambda}^{n-1}$ and $S_{\eta}^{n-1}$ are the unit sphere in $\mathbb{R}_{\lambda}^n$ and $\mathbb{R}_{\eta}^n$ respectively.
 
From \eqref{4.56} it follows that $\Sigma_{(\gamma_k)}\cap S_{\lambda}^{n-1}$ is a compact set and so is $T_{(\gamma_k)}.$ Since $\widetilde{ d^{2} f_{k}(\lambda, \eta)}$ is negative and continuous on $T_{(\gamma_k)},$ there exists $\gamma^{(k)}>0$ such that $\widetilde{d^2f_k(\lambda,\eta)}\le -\gamma^{(k)} $ on $T_{(\gamma_{k}) }$ and the theorem is proved.
		\end{proof}
\section{The $d-$concavity}
\subsection{Completing estimates for $d^2F_k(R,M)$}
Suppose $0<\delta<1, \mu >0, $ $0< \gamma_k <1.$ We denote 
\begin{equation}\label{5.1}
	D_{\delta,\mu,\gamma_k}=\left\{R=\omega+\beta\in D_{\delta,\mu}; \lambda(\omega)\in \Sigma_{(\gamma_k)} \right\}.
\end{equation}
It is clear that $D_{\delta,\mu, \gamma_{k}}$ is a convex and unbounded set.
\begin{proposition}\label{md5.1}
	Suppose $2\le k\le n-1.$ There exists $0<\delta_0<1,$ that depends only on $k, n, \gamma_k$ and there exists $C_{10}	>0$ that  depends 
 only on $ k, n, \gamma_k, \delta_{0}$ and does not depend on $\mu$ such that for all $R \in D_{\delta, \mu, \gamma_k}$ and all $P,$ $P^T=P$ the following estimate holds  for all $\delta: 0< \delta \le \delta_0$
 \begin{equation}\label{5.2}
 	d^2F_k(R,P)\le -C_{10} | \tilde{\tilde{P}}|^2.
 \end{equation}
\end{proposition}
\begin{proof}
From \eqref{2.17}, 	\eqref{3.26} we have
\begin{equation}\label{5.3}
	d^2F_k(R,P)=d^2F_k(\widetilde{R},\widetilde{P})\le d^2F_k(D,\widetilde{P})+C_4\delta^2 | \tilde{\tilde{P}}|^2.
\end{equation}
From 	\eqref{3.48} we also have 
\begin{equation}\label{5.4}
	\begin{split}
	d^2F_k(D,\widetilde{P})&=d^2F_k(D,\widetilde{U})	-\frac{1}{S_k(D)} \sum_{i_1, \cdots, i_k} G_{i_1 \cdots i_k}(D) \sum_{p,q=1\atop p\ne q}^k  | \tilde{\tilde{P}}_{i_pi_q}|^2\\
	&= d^2F_k(D,\widetilde{U}) -\sum_{\ell, m=1\atop \ell\ne m}^n \left[ \sum_{i_1,\cdots,i_k\atop \ell, m\in [i_k] } \gamma_{i_1\cdots  i_k}\right]| \tilde{\tilde{P}}_{\ell m}|^2,
	\end{split}
\end{equation} 
where $ \gamma_{i_1\cdots  i_k}=\frac{\lambda_{i_1} \cdots \lambda_{i_k}}{\sigma_{k}(\lambda)}, \lambda=\lambda(\omega)=(\lambda_1,\cdots, \lambda_n),$  $\lambda(\omega)\in  \Sigma_{(\gamma_{k})},$ $\omega=C^{-1}DC,$ $C^T=C^{-1},$ $D={\rm diag} (\lambda_1,\cdots, \lambda_n), $ $\widetilde{P}=CPC^{-1}, $ $ \tilde{\tilde{P}}=D^{-\frac{1}{2}} \widetilde{P}  D^{-\frac{1}{2}},$ $\widetilde{P}=[\widetilde{P}_{ij}]_{n\times n},$ $\widetilde{U}={\rm diag} (\widetilde{P}_{11},\cdots, \widetilde{P}_{nn}).$

From \eqref{4.58T}, \eqref{4.8} we get 
\begin{equation}\label{5.5}
	\begin{split}
	d^2F_k(D,\widetilde{U})&=	d^2f_k(\lambda,\xi)=d^2f_k(\lambda,\lambda\frac{\xi}{\lambda})=\widetilde{ d^{2} f_{k}(\lambda, \eta)}\\
	&\le -\gamma^{(k)}|\eta|^2=-\gamma^{(k)}\sum_{i=1}^{n} \left( \frac{\xi_i}{\lambda_i}\right)^2=- \gamma^{(k)}\sum_{i=1}^{n} \left( \frac{\widetilde{P}_{ii}}{\lambda_i}\right)^2=-\gamma^{(k)}\sum_{i=1}^{n} |\tilde{\tilde{P}}|^2,
	\end{split}
\end{equation}
where $\eta=\frac{1}{\lambda}\xi =\left( \frac{\xi_1}{\lambda_1}, \frac{\xi_2}{\lambda_2},\cdots  \frac{\xi_n}{\lambda_n}\right),$ $\xi =(\xi_1,\cdots, \xi_n)=\left( \widetilde{P}_{11},\cdots, \widetilde{P}_{nn}\right).$
\end{proof}
To estimate the second term on the right-hand side of \eqref{5.4} we need the following 
\begin{lemma}\label{bd5.1}
	Suppose $2\le k\le n-1$ and $\lambda(D)=(\lambda_1,\cdots, \lambda_n)\in \Sigma_{(\gamma_k)}.$ Then  the following inequality holds 
\begin{equation}\label{5.6}
 	\sum_{   i_1,\cdots,  i_k \atop   \ell, m\in \{i\}_k, \ell\ne m}\gamma_{i_1\cdots  i_k}  \ge \mu^{(k)} =\frac{(k-1)k}{(n-1)n} (\gamma_k)^k.
\end{equation}	
	\end{lemma}
\begin{proof}
 We have 
 $$\sum_{   i_1,\cdots,  i_k \atop   \ell, m\in \{i\}_k, \ell\ne m}\lambda_{i_1}\cdots \lambda_{i_k} \ge \binom{n-2}{k-2}\lambda_n^k,\quad \sigma_{k}(\lambda)\le \binom{n}{k}\lambda_1^k.$$
Therefore, 
$$\sum_{   i_1,\cdots,  i_k \atop   \ell, m\in \{i\}_k, \ell\ne m}\gamma_{i_1\cdots i_k} \ge \frac{(k-1)k}{(n-1)n}(\gamma_k)^k=\mu^{(k)}. $$	
	
Continuation of proof for Proposition \ref{md5.1}: It is clear that from \eqref{5.4}-\eqref{5.6} we can choose $C_{11}=\min(\gamma^{(k)},\mu^{(k)})$ to get $d^2F_k(D,\widetilde{P})\le -C_{11} |\tilde{\tilde{P}}|^2$ and 
\begin{equation*} 
	d^2F_k(D,P) \le - C_{11} |\tilde{\tilde{P}}|^2+C_4\delta^2 |\tilde{\tilde{P}}|^2=-\frac{1}{2}C_{11} |\tilde{\tilde{P}}|^2+ (C_4\delta^2-\frac{1}{2}C_{11}) |\tilde{\tilde{P}}|^2,
\end{equation*}
from which we choose $\delta_0=\min\left(\frac{1}{2},\sqrt{\frac{C_{11}}{2C_4}}\right), C_{10}=\frac{1}{2}C_{11}.$
\end{proof}
Now we formulate one of  main results of the paper.
\begin{theorem}\label{dl5.3T}
	Suppose $ 2\le k\le n-1$ and $0<\delta_0<1$ as was described in Proposition \ref{md5.1}. Then there exists  $0<\delta_1\le \delta_0,$ $C_{12}>0$ that depends on $k, n, \gamma_k, \delta_{0}$ and does not depend on $\mu$ such that for all $R\in D_{\delta,\mu,\gamma_k}$ and all $M=P+Q, P^T=P, $ $Q^T=-Q$ the following estimate holds for $\delta:$ $0<\delta\le \delta_1$
	\begin{equation}\label{5.10}
	d^2F_k(R,M)	\le C_{12} \frac{|Q|^2}{\lambda^2_{\min}(\omega)}.
	\end{equation}
	\end{theorem}
\begin{proof}
 From \eqref{2.17}, \eqref{3.47},  we have
 \begin{equation*}\label{5.11}
\begin{split}
 	d^2F_k(R,M)=	d^2F_k(\widetilde{R},\widetilde{M})  	&\le 	d^2F_k(D,\widetilde{P})+(C_4+1)\delta^2 |\tilde{\tilde{P}}|^2+ (C_9^2+C_8)|\tilde{\tilde{Q}}|^2\\
 &\le - C_{11} |\tilde{\tilde{P}}|^2 +(C_4+1)\delta^2 |\tilde{\tilde{P}}|^2+ (C_9^2+C_8)|\tilde{\tilde{Q}}|^2.
 \end{split} 
 \end{equation*}
	Then we choose $\delta_1$ so that  $\delta_1=\min\left(\delta_{0},\sqrt{\frac{C_{11}}{C_4+1}}\right) $
	and $C_{12}=C^2_9+C_8.$ We note that 
$$|\tilde{\tilde{Q}}|=|D^{-\frac{1}{2}} \widetilde{Q}  D^{-\frac{1}{2}}|\le \frac{ |\tilde{Q}|}{\lambda_{\min}^2(D^{\frac{1}{2}})}=\frac{|Q|}{\lambda_{\min}(\omega)}.$$
	\end{proof}
\subsection{The $d-$concavity}
As consequence of Theorem \ref{dl5.3T}, we formulate below two variants of $d-$concavity of the function $F_k(R)$ on the set $D_{\delta,\mu,\gamma_k}.$
\begin{theorem}\label{dl5.2}
	Suppose $2\le k\le n-1$ and  $R^{(0)}, R^{(1)}\in D_{\delta,\mu,\gamma_k},$ where $R^{(0)}=\omega^{(0)}+\beta^{(0)},$ $R^{(1)}=\omega^{(1)}+\beta^{(1)}.$ Then we have
\begin{equation}\label{5.12}
		F_{k}\left(R^{(1)}\right)-F_{k}\left(R^{(0)}\right)\le \sum_{i,j=1}^n  \frac{\partial F_{k}\left(R^{(0)}\right)}{\partial R_{ij} }\left(R^{(1)}_{ij}-R_{ij}^{(0)}\right)	+C_{12} \frac{| \beta^{(1)}-\beta^{(0)}|^{2}}{\lambda^2_{\min}\left(\omega^{(\tau)}\right),}
\end{equation}
where $\omega^{(\tau)}=(1-\tau)\omega^{(0)}+\tau\omega^{(1)},\ 0<\tau<1.$
		\end{theorem}
\begin{proof}
We consider the following function on $[0,1]:$
\begin{equation*}\label{5.13}
	g(t)=F_{k}\left(R^{(t)}\right),
\end{equation*}
where $R^{(t)}=(1-t) R^{(0)}+tR^{(1)}.$ We have 
\begin{equation*}\label{5.15}
	\begin{split}
	g'(t)&=\sum_{i, j=1}^{n} \frac{\partial F_{k}\left(R^{(t)}\right)}{\partial R_{ij}}\left(R_{i j}^{(1)}-R_{i j}^{(0)}\right),\\
	g"(t)&=\sum_{i,j, \ell, m} \frac{\partial^{2} F_{k}(R(t))}{\partial R_{ij}  \partial R_{\ell m}} \left(R_{i j}^{(1)}-R_{\ell m}^{(0)}\right)\left(R_{ij}^{(1)}-R_{\ell m}^{(0)}\right),
	\end{split}
\end{equation*}
\begin{equation}\label{5.16}
	g(1)-g(0)=g'(0)+g"(\tau),\quad 0< \tau <1.
\end{equation}
It is obvious that 
\begin{equation}\label{5.17}
g(1)=F\left(R^{(1)}\right),\quad 	g(0)=F\left(R^{(0)}\right).
\end{equation}
From \eqref{5.10} we have
\begin{equation}\label{5.18}
	g''(\tau)=d^{2} F_{k}\left(R^{(\tau)}, R^{(1)}-R^{(0)}\right) \le C_{12} \frac{\left|\beta^{(1)}-\beta^{(0)}\right|^{2}}{\lambda_{\min }^{2}\left(\omega^{(\tau)}\right)}.
\end{equation}
From \eqref{5.16}-\eqref{5.18} we get \eqref{5.12}.
	\end{proof}
\begin{theorem}\label{dl5.3}
Under the assumptions of the 	Theorem \ref{dl5.2} we have
\begin{equation}\label{5.19}
	F_{k}\left(R^{(1)}\right)-F_{k}\left(R^{(0)}\right)\le \sum_{i, j=1}^{n} \frac{\partial F_{k}\left(R^{(0)}\right)}{\partial R_{i j}} \left(R_{1}^{(1)}-R_{1 j}^{(0)}\right)+d,
\end{equation}
where $d=4nC_{12}\delta^2.$ This means that the function $F_{k}(R)=\log \left(S_{k}(R)\right)$ is $d-$concave on the set $D_{\delta,\mu,\gamma_k}.$
	\end{theorem}
\begin{proof}
  We use the inequality \eqref{5.12}. We have 
  $$|\beta^{(1)}- \beta^{(0)} | \le \sqrt{n} \|\beta^{(1)}- \beta^{(0)} \|\le 2\sqrt{n} \mu, $$
  $$\mu \le \delta\lambda_{\min }(\omega^{(\tau)}),$$
  from which it follows that we can choose $d=4nC_{12}\delta^2$ in \eqref{5.19}. 
	\end{proof}

 \section*{Acknowledgment } 
This research was funded by Hanoi Pedagogical University 2 Foundation for Sciences and Technology Development via grant number: C2020-SP2-05.

\end{document}